\documentclass[11pt,reqno]{amsart}

\usepackage[dvipsnames]{xcolor}
\usepackage{amsmath}
\usepackage{amsthm}
\usepackage{graphicx}
\usepackage{upgreek}
\usepackage{hyperref}
\usepackage{mathrsfs}
\usepackage{bm,bbm}
\usepackage{amsfonts}
\usepackage{amssymb}
\usepackage{csquotes}
\usepackage{stmaryrd}
\usepackage{ulem}
\usepackage{tikz-cd}
\usepackage{cancel}
\usepackage{siunitx}
\usepackage{algorithm}
\usepackage{algpseudocode}
\usepackage{caption}
\usepackage{subcaption}
\usepackage{booktabs}
\usepackage{enumerate}
\usepackage{soul}
\usepackage[margin=1in]{geometry}
\numberwithin{equation}{section}

\theoremstyle{plain}
\newtheorem{theorem}{Theorem}[section]
\newtheorem{definition}[theorem]{Definition}
\newtheorem{proposition}[theorem]{Proposition}

\newtheorem{remark}[theorem]{Remark}
\newtheorem{example}[theorem]{Example}
\newtheorem{corollary}[theorem]{Corollary}

\newcommand{\thistheoremname}{}
\newtheorem{genericthm}[theorem]{\thistheoremname}

\newtheorem*{genericthm*}{\thistheoremname}
\newenvironment{assumption*}[1]
  {\renewcommand{\thistheoremname}{#1}%
   \begin{genericthm*}}
  {\end{genericthm*}}

\newcommand{\R}{\mathbb{R}}

\renewcommand{\P}{\mathbb{P}}

\newcommand{\sC}{\mathcal{C}}

\newcommand{\sL}{\mathcal{L}}
\newcommand{\N}{\mathbb{N}}

\newcommand{\E}{\mathbb{E}}

\newcommand{\lbrs}[1]{\llbracket #1 \rrbracket}

\newcommand{\intr}{\text{Int}}

\title{Convergence of the Deep Galerkin Method for Finite State Mean Field Control Problems}

\author{William Hofgard}
\address{Department of Mathematics, University of California, Berkeley, CA}
\email{hofgard@berkeley.edu}

\author{Jingruo Sun}
\address{Department of Operations Research and Financial Engineering, Princeton University, Princeton, NJ}
\email{jingruo@princeton.edu}

\author{Asaf Cohen}
\address{Department of Mathematics, University of Michigan, Ann Arbor, MI}
\email{shloshim@gmail.com}
\urladdr{https://sites.google.com/site/asafcohentau/}

\date{\today}

\thanks{This is the final version of the paper, to appear in the SIAM Journal on Mathematics of Data Science.}

\begin{document}
\maketitle
\begin{abstract}
We establish the convergence of the deep Galerkin method (DGM), a deep learning-based scheme for solving high-dimensional nonlinear PDEs, for Hamilton-Jacobi-Bellman (HJB) equations that arise from the study of mean field control problems (MFCPs). Based on a recent characterization of the value function of the MFCP as the unique viscosity solution of an HJB equation on the simplex, we establish both an existence and convergence result for the DGM. First, we show that the loss functional of the DGM can be made arbitrarily small given that the value function of the MFCP possesses sufficient regularity. Then, we show that if the loss functional of the DGM converges to zero, the corresponding neural network approximators must converge uniformly to the true value function on the simplex. We also provide numerical experiments demonstrating the DGM's ability to generalize to high-dimensional HJB equations.
\end{abstract}

\noindent{\bf AMS Classification:} 
91A07,  	%Games with infinitely many players
35Q89,  	%PDEs in connection with mean field game theory
68T07,  	%Artificial neural networks and deep learning
49L12,  	%Hamilton-Jacobi equations in optimal control and differential games
49N10, 	%Linear-quadratic optimal control problems
35A35,  	%Theoretical approximation in context of PDEs
60J27 %Continuous-time Markov processes on discrete state spaces
%35M99, % PDEs 	None of the above, but in this section   OR MAYBE :35M10  	Equations of mixed type
%93E20, %	Optimal stochastic control
%60F05, % 	Central limit and other weak theorems
%91A15 % 	Stochastic games

%\tableofcontents

\section{Introduction}
We consider the convergence of a deep learning-based numerical method for solving Hamilton--Jacobi--Bellman (HJB) equations that arise from the study of mean field control problems (MFCPs). The specific class of HJB equations that we consider are first-order, nonlinear PDEs on the $(d-1)$--dimensional simplex $S_d$, of the following form:
\begin{align}
\label{eq:mfcp-hjb}
\begin{split}
    &-\partial_t V(t, m) + \sum_{i \in \lbrs{d}} m_i H^i(t, m, D^i V(t,m)) = 0,\\
    &V(T, m) = \sum_{i \in \lbrs{d}} m_i g^i(m).
\end{split}
\end{align}
Above, $t \in [0, T]$ for some $T > 0$, $m \in S_d$ lies in the simplex, $\lbrs{d} := \{1, \ldots, d\}$, $V : [0,T] \times S_d \to \R$ is the {\it value function} (optimal cost) of an associated MFCP, $H^i$ is a Hamiltonian (to be introduced below in the context of a MFCPs), and $g^i : S_d \to \R$ is a terminal cost defined on the simplex. The operator $D^i$ denotes a a derivative along the simplex.
Namely, $$D^i V(t,m)=(\partial_{m_j - m_i} V(t,m))_{j \in \lbrs{d}} = (\nabla_m V(t, m) \cdot (e_j - e_i))_{j \in \lbrs{d}},$$
where $e_i$ denotes the $i$th standard basis vector in $\R^d$. Note that on the simplex, directional derivatives are only permitted in the directions $e_j - e_i$. Consequently, we can write Equation~\eqref{eq:mfcp-hjb} in standard PDE notation:
\begin{align}
\label{eq:mfcp-hjb-pde}
\begin{split}
    &\partial_t V(t, m) = F(t, m, \nabla_m V(t,m)) := \sum_{i \in \lbrs{d}} m_i H^i(t, m, D^i V(t,m)),\\
    &V(T, m) = G(m) := \sum_{i \in \lbrs{d}} m_i g^i(m).
\end{split}
\end{align}
We focus on solving Equation~\eqref{eq:mfcp-hjb} numerically due to the connection between the MFCP and the $N$-agent optimization problem. 

Our numerical approach is primarily based on the deep Galerkin method (DGM), as first introduced by \cite{dgm}. In particular, we solve the above class of PDEs using a deep learning-based algorithm motivated by the classical finite element method. However, instead of finding a basis of functions to approximate solutions as in the classical method, we utilize a sufficiently rich class of neural networks to approximate solutions to Equation~\eqref{eq:mfcp-hjb}. In addition to proving the convergence of our numerical scheme under the assumptions outlined in Section~\ref{sec:hjb-mfcp-setup}, we provide numerical experiments that demonstrate the validity of this approach. First, however, we introduce the MFCP model in Section~\ref{subsec:mfcp}, describe the current field of solving high-dimensional PDE using deep learning in Section~\ref{subsec:num-sol-hpde}, and outline our main results in Section~\ref{subsec:main-results}.

\subsection{Mean Field Control Problems}
\label{subsec:mfcp}
\subsubsection{Literature review}
MFCPs are limiting models for cooperative games with a finite, but large, number of interacting agents attempting to minimize a common cost. Namely, the study of MFCPs is motivated by situations in which a large number of agents aim to achieve a common goal. For instance, applications arise in the control of autonomous vehicles and drones \cite{drones7070476, mot-pha2022CONTROL}, efficient real-time streaming between devices \cite{7218597-2015}, decentralized and centralized crowd motion \cite{crowd-motion-1, crowd-motion-2}, analysis of wireless local-area networks (WLANs) \cite{duffy2010meanCONTROL, bor-sun2012CONTROL}, among many other related applications. The inception of limiting models for many-agent games can be traced back to Aumann's seminal work in 1964 \cite{aumann1964markets} and Schmeidler's influential contribution in 1973 \cite{schmeidler1973equilibrium}.  These works dealt with one-shot games. In the context of stochastic dynamical games, mean field game (MFG) theory was initially developed Lasry and Lions \cite{lasry-lions-mfg} and independently by Huang, Caines, and Malhamé \cite{huang-caines-malhame-mfg}. In our discourse, we distinguish between MFCPs and MFGs, which are frameworks for approximating, respectively, cooperative and non-cooperative games involving numerous players. MFCPs represent optimal control problems for McKean--Vlasov dynamics, where the value function (optimal cost) is often delineated by a HJB equation. Conversely, in MFGs, the Nash equilibrium is approximated through a fixed point of a best response map, and the value function (i.e., the cost under the proper equilibrium concept in the MFG) is characterized by what is called the master equation. 

Recent work by Pham and Wei derived a dynamic programming principle and HJB equation for the MFCP by framing it as a deterministic control problem of a Fokker--Planck equation; see \cite{Pham2017, Pham2018}. Even more recently, work by Cecchin \cite{mfcp} establishes important existence and uniqueness results for the HJB equation corresponding to the finite-state MFCP we consider here. For a comprehensive overview of MFG theory, MFCPs, and their relationship to stochastic games and control problems, see \cite{Bensoussan2013, Carmona2018_I, Carmona2018_II}.

Given current advances in deep learning, there has been recent progress towards solving both MFGs and MFCPs numerically using deep learning \cite{ml-mfcp, deep-learning-mfcp}. Since 2020, several review articles have again brought up the potential for solving MFGs and MFCPs with deep learning methods, including the DGM and deep BSDE; see \cite{mfc-mfcp-review, scalable-rl-mfg}. For MFCPs with \textit{continuous} state spaces, for example, \cite{Carmona2021ergodic} and \cite{carmona2021finitehorizon} consider the convergence of machine learning-based algorithms in the ergodic and finite horizon cases respectively. For finite state MFGs, \cite{cohen2024deep} provides an analysis of convergence of both DGM-type and deep BSDE approaches. Similarly, \cite{actor-critic-bsde} provides a broad sampling of applications of deep BSDE-like algorithms to common problems in stochastic control without providing a convergence guarantee for their proposed method. Building on the theoretical results from \cite{mfcp}, we aim to provide a novel pathway towards solving finite-state MFCPs using deep learning, complete with a robust convergence guarantee.

\subsubsection{MFCP model} The $N$-agent optimization problem, which limits to the MFCP as $N$ approaches infinity, is a continuous-time, finite state optimization problem with $N$ agents whose state process is $(\mathbf{X}_s = (X^1_s, \ldots, X^N_s))_{s\in[t_0,T]}$, where each coordinate takes values in the state space $\lbrs{d} = \{1, \ldots, d\}$. Above, $[t_0,T]$ is the deterministic time interval of the problem. The agents' state process $(\mathbf{X}_t)_{t\in[t_0,T]}$ follow the dynamics given by the continuous-time Markov chain, where each of its coordinates follows 
\begin{align}
\label{eq:n-agent}
    \P(X_{s+h}^k = j \mid \mathbf{X}_s = \mathbf{x}) = \beta^k(s, \mathbf{x};j)h + o(h)\qquad\text{as \qquad $h\to0^+$},\qquad \forall j\in\lbrs{d},\; j\ne x^k.
\end{align}
The agents choose feedback controls $(\beta^1, \ldots, \beta^N)$ in order to minimize the common cost
\begin{align}
    J^N(t_0, m^N_{t_0}, \beta^1, \ldots, \beta^N) := \frac{1}{N} \sum_{k=1}^N\E\left[\int_{t_0}^T f(s, X_s^k, \beta^k(s, \mathbf{X}_s;\cdot), \mu_s^N) ds + g(X^k_T, \mu_T^N)\right] \label{eq:n-agent-cost},
\end{align}
where $\mu^N_s$ is the empirical distribution of the $N$ agents at time $s$, satisfying the predetermined initial condition $\mu^N_{t_0}=m^N_{t_0}$, $f$ denotes the running cost of the stochastic control problem while $g$ denotes the terminal cost. The value of $f(s,i,a,m)$ is independent of $a_i$ for each $i \in \lbrs{d}$ because there is no transition rate between any state $i \in \lbrs{d}$ and itself. 

In practice, finding the optimal control in $N$-agent stochastic control problems quickly becomes intractable as $N$ grows. Due to the connection between the $N$-agent problem and the MFCP established in \cite{mfcp}, for instance, an efficient method for computing the value function and optimal control in the MFCP therefore provides a path for solving difficult $N$-agent stochastic control problems. Before presenting the MFCP, we provide motivation for its structure. As $N$ approaches infinity, the $N$-agent optimization problem resembles an optimization problem involving a single {\it representative} agent whose dynamics follow a controlled continuous-time Markov chain of the form:
\begin{align}
\label{eq:mfcp-process}
    \P(X_{s+h} = j \mid X_s = i) = \alpha_{i,j}(s) h + o(h)\qquad\text{as}\qquad h\to0^+,
\end{align}
where the feedback control of the representative agent $\alpha_{i,j}(s)$ is a measurable, Markovian function, standing for the transition rate of the above process from state $i$ to $j$, both in $\lbrs{d}$.
Throughout this paper, we consider only rate matrices satisfying $\alpha_{i,j}(s) \in [0, M]$ for some constant $M > 0$ and all $i\ne j, i,j \in\lbrs{d}$, $s\in[t_0,T]$, and we refer to them as {\it admissible controls}. We denote the set of admissible controls by $\mathcal{A}_{t_0}$. Also, $\alpha_i$ refers to the $i$th row of the transition matrix $\alpha$. 

At this point, the reader may wonder about the gap between the richness of the feedback controls allowed in the $N$-agent problem versus those allowed the MFCP. Specifically, the feedback controls in the $N$-agent problem may depend simultaneously on the states of all $N$ agents, but the feedback controls in the MFCP only depend on the state of the representative agent. We posit that, despite this apparent disparity, there is not a significant difference between the classes of controls considered in the two models. Indeed, in the $N$-agent model, all agents are symmetric. Given the symmetry among agents, a logical step involves substituting their transition rates with functions of the form $\bar \beta^k(s,X^k_s,\mu^N_s;\cdot)$ that depend on the empirical distribution of the $N$ agents instead of all of their states. Guided by the propagation of chaos result from \cite{mfcp}, we anticipate the convergence of the empirical distributions $(\mu^N_t)_{t\in [t_0,T]}$ toward a deterministic flow of measures, approximating the distributions of each coordinate in $\mathbf{X}$. This deterministic framework implies that, in the limit, we can anticipate the control's dependence on both the current time $t$ and the present state $i$, exactly as is the case in the MFCP.

As in the $N$-agent case, the representative agent aims to minimize the cost
\begin{align}
    J(t_0, m_{t_0}, \alpha) := \E\left[\int_{t_0}^T f(s, X_s, \alpha_{X_s}(s), \text{Law}(X_s))ds + g(X_T, \text{Law}(X_T))\right] \label{eq:mfcp-cost},
\end{align}
where the initial distribution of the process $X$: $\text{Law}(X_{t_0}) = m_{t_0}\in S_d$ is predetermined. This problem can be alternatively described as a deterministic control problem in terms of the Fokker--Planck equation for $(\mu_s := \text{Law}(X_s))_{s\in[t_0,T]}$. In particular, for a given admissible control $\alpha$, $\mu$ uniquely solves the ordinary differential equation (ODE):
\begin{align}
\label{eq:forward-eq}
\begin{split}
    &\frac{d}{dt}\mu_t^i = \sum_{j \in \lbrs{d}} \left(\mu_t^j \alpha_{j,i}(s) - \mu_t^i \alpha_{i,j}(s)\right),\qquad i\in\lbrs{d}, \quad t\in[t_0,T], \\
    &\mu_{t_0} = m_{t_0},
\end{split}
\end{align}
and the cost functional for the deterministic control problem becomes
\begin{align}
    \label{eq:det-mfcp-cost}
    J(t_0, m_{t_0},\alpha) := \int_{t_0}^{T} \sum_{i \in \lbrs{d}} f(s, i, \alpha_{i}(s), \mu_s) \mu_s^i ds + \sum_{i \in \lbrs{d}} g^i(\mu_T) \mu_T^i.
\end{align}
Allowing for the initial time of the problem $t_0$ to be in $[0,T]$, the {\it value function} $V:[0,T]\times S_d\to\R$ of the deterministic control problem is then defined by
\begin{align}
    \label{eq:value-function}
    V(t, m) = \inf_{\alpha \in \mathcal{A}_{t}} J(t, \alpha, m).
\end{align}
By standard optimal control arguments via an appropriate dynamic programming principle, this gives rise to the HJB equation \eqref{eq:mfcp-hjb}, where the {\it PDE-Hamiltonian} $H : [0, T] \times S_d \times \R^d \to \R$ is defined by $H(t, m, z) = \sum_{i \in \lbrs{d}} m_i H^i(t, m, z)$, with
\begin{align}
\label{eq:hamiltonian}
    H^i(t,m,z) := \sup_{a \in [0, M]^d, \;k\in\lbrs{d}\setminus\{i\}} \Big(-\sum_{k\in\lbrs{d}\setminus\{i\}}a_k z_k - f(t,i,a,m)\Big).
\end{align}
Note that we refer to both $H^i$ and $H$ as Hamiltonians. To avoid confusion, we occasionally employ the term PDE-Hamiltonian for $H$ in particular. Under the assumptions presented in Section~\ref{subsec:assumptions} below, the above equation has a unique maximizer \cite{unique-optimal-control-mfg}, which we denote by $a^*(t, i, m, z)$ for each $i \in \lbrs{d}$. In turn, the \textit{unique} optimal control for the MFCP is given by
\begin{align}
    \label{eq:optimal-control-mfcp}
    \alpha^*_{i,j}(t,m) = a_j^*(t, i,  m, D^i V(t, m)),
\end{align}
where the directional derivative of the value function $D^i V(t, m)$ on the simplex $S^d$ is as defined in Equation~\eqref{eq:mfcp-hjb}. Additionally, the dynamics of the process $\mu$ under the optimal control satisfy the dynamics in Equation~\eqref{eq:forward-eq} upon replacing the control therein with the optimal control from Equation~\eqref{eq:optimal-control-mfcp} above \cite{mfcp}.

Importantly, recent work rigorously establishes the connection between this $N$-agent stochastic control problem and the corresponding MFCP. In particular, Cecchin \cite{mfcp} formulates the MFCP in terms of the deterministic optimal control of a Fokker--Planck equation as in Equations~\eqref{eq:forward-eq} and \eqref{eq:det-mfcp-cost} above. The resulting deterministic control problem yields a dynamic programming principle and the HJB equation \eqref{eq:mfcp-hjb} solved by the value function of the MFCP.  By formulating the $N$-agent optimization problem in a similar manner, one obtains two HJB equations: one for the $N$-agent optimization problem and one for the MFCP. Both HJB equations have unique viscosity solutions, as shown in \cite{mfcp}. In turn, \cite{mfcp} establishes an explicit rate of convergence between $V^N$, the value function for the $N$-agent problem, and $V$, the value function for the MFCP. In particular, a convergence rate of $1/\sqrt{N}$ is established. In a similar vein, Kolokoltsov \cite{kolokoltsov-conv} obtains a convergence rate of $1/N$ under additional regularity assumptions. Notably, \cite{mfcp} places no regularity assumptions on the value function $V$ (other than Lipschitz continuity), and in the most general case, does not require convexity of the running or terminal costs of the MFCP. Given the relationship between the $N$-agent problem and the MFCP, by approximately solving the HJB equation for the MFCP in this paper, we essentially provide a method for approximately solving the $N$-agent problem.

\subsection{Numerical Solutions to High-Dimensional PDE}
\label{subsec:num-sol-hpde}
We aim to construct a numerical scheme for efficiently solving the HJB equation associated with the MFCP. In particular, as the dimension $d$ increases, the so-called ``curse of dimensionality'' prevents standard numerical schemes (e.g., Monte Carlo methods, mesh-based algorithms, etc.) from solving the HJB equation in a tractable manner. However, recent advancements in deep learning present promising options for solving high-dimensional, non-linear PDE such as the HJB equation in Equation~\eqref{eq:mfcp-hjb}.

Two leading methods have been presented for parabolic PDEs that resemble the HJBs equation: the deep Galerkin method (DGM) and deep backward stochastic differential equations (BSDE). The primary focus of this paper is the DGM. This approach, first introduced in \cite{dgm}, models itself after classical finite element methods for solving low-dimensional PDEs. However, the DGM is a mesh-free method; instead of creating basis functions that approximate the solution to a PDE from a mesh, the DGM utilizes neural network approximations that only depend on the parameters, the architecture, and the activation function between layers of the network. The loss functional of the DGM with $L^2$-loss attempts to minimize the $L^2$-error of both the PDE and the terminal condition during training, ultimately learning the parameters that best approximate the solution of the PDE. In \cite{dgm}, the authors introduce the DGM, illustrate its ability to numerically solve high-dimensional nonlinear PDEs, and provide a convergence guarantee for second-order nonlinear parabolic PDEs.

The second popular deep learning-based method for solving high-dimensional PDEs, introduced in \cite{deep-bsde-orig} and expanded upon in \cite{nonlinear-feynman-kac, db-schemes}, exploits the connection between nonlinear parabolic PDEs and backward SDEs can be exploited via a so-called non-linear Feynman--Kac formula. In turn, the resulting backward SDE can be solved numerically by recursively using a sequence of neural networks to solve the SDE along a specified discretization of the time interval in question, starting with the terminal condition of the original PDE. Although this method, often referred to as deep BSDE, is likely applicable to our context, we defer further consideration of deep backward schemes to future work.

More recently, several works have investigated the problem of convergence when solving PDEs with deep learning methods. For instance, \cite{jiang2023globalconvergencedeepgalerkin} show that an extension of the common $Q$-learning algorithm from reinforcement learning results in neural network approximators that solve a broad class of elliptic PDEs in the infinite-width limit. Using a similar neural tangent kernel (NTK) analysis, \cite{jiang2023globalconvergencedeepgalerkin} shows that DGM produces approximations that converge to the true solutions of second-order linear PDEs. Finally, \cite{cohen2024deep} studies convergence of both DGM-type methods and deep BSDE for first-order nonlinear PDEs arising from finite state MFGs.

\subsection{Main Results}
\label{subsec:main-results}
We show that the convergence guarantee provided in \cite{dgm} extends to the (first-order) HJB equation \eqref{eq:mfcp-hjb}, relying on the theory of viscosity solutions for HJB equations to obtain the desired convergence. Importantly, our class of HJB equations does not fall into the class of second-order nonlinear parabolic PDEs considered in \cite{dgm}, so we provide a different proof technique. In order to leverage the theory of viscosity solutions in our proof technique, we modify the loss functional used in \cite{dgm} to approximate the uniform norm rather than the $L^2$-norm. Namely, we choose an $L^\infty$-loss. Indeed, viscosity solutions are often referred to as $L^\infty$-weak solutions, as they respect convergence with respect to the uniform norm \cite{crandall-evans-lions}. The work by \cite{dgm} establishes the uniform convergence of neural network approximators to the actual solution a family of partial differential equation. However, this result is obtained by introducing additional assumptions of uniform boundedness and equicontinuity on the neural networks, whch we bypass using the $L^\infty$-loss. With this reformulation of the DGM loss, we first prove an existence result by showing in Theorem~\ref{th:unif-error-conv} that a sequence of neural networks taking the loss functional of the DGM to zero exists. We then prove in Theorem~\ref{th:unif-value-conv} that for such a sequence of networks, the neural network approximators converge uniformly to the true value function for the MFCP on $[0, T] \times S_d$. 

To obtain this result, we utilize a standard argument from the theory of viscosity solutions, relying on an appropriate definition of viscosity solutions to Equation~\eqref{eq:mfcp-hjb} on the simplex and a corresponding comparison principle. By establishing that a sequence of neural network approximators are viscosity solutions to a sequence of perturbed HJB equations and applying the comparison principle to two suitable upper and lower limits of the neural network approximators, we can conclude that the sequence of neural network approximators taking the DGM loss functional to zero must converge uniformly to the classical solution to Equation~\eqref{eq:mfcp-hjb}. The main steps in our proof are outlined as follows:
\begin{itemize}
    \item First, in Theorem~\ref{th:unif-error-conv}, we show that there exists a neural network $\varphi$ that makes the DGM loss (defined in Section~\ref{sec:dgm}) arbitrarily small. 
    \item From this, in Corollary~\ref{corr:convergence-to-v}, we establish that there exists a sequence of neural network approximators $\{\varphi^n\}_{n \in \N}$ such that $\varphi^n \to V$ uniformly on $[0,T]\times S_d$ as $n \to \infty$, where $V$ is the unique classical solution to Equation~\eqref{eq:mfcp-hjb}.
    \item Then, we handle the technical details of working with a PDE on the $(d-1)$--dimensional simplex $S_d$. In particular, the interior of $S_d$ is empty, so we instead work on $\widehat{S}_d$, the projection of the simplex onto $\R^{d-1}$. These technical details are discussed in Proposition~\ref{prop:eps-n-convergence}, and we work on $\widehat{S}_d$ for the remainder of the proof in order to fully leverage the power of viscosity solutions. Throughout, we use the {\it hat} notation to denote mathematical objects that are defined on $\widehat{S}_d$ rather than on $S_d$. Specifically, $\widehat{V}$ is the version of $V$, defined on $\widehat{S}_d$.
    \item Finally, we arrive at our main result in Theorem~\ref{th:unif-value-conv}, which establishes that if $\{\widehat{\varphi}^n\}_{n \in \N}$ is a sequence of neural network approximators taking the DGM loss to zero (which we know exists from Corollary~\ref{corr:convergence-to-v}), then $\widehat{\varphi}^n \to \widehat{V}$ uniformly on $[0, T] \times \widehat{S}_d$. 
    \item The proof of Theorem~\ref{th:unif-value-conv} relies on a comparison principle for viscosity solutions, as presented in \cite[Theorem 3.3]{visc-users}. In particular, we construct a viscosity subsolution $\overline{V}$ and a viscosity supersolution $\underline{V}$ on $[0, T] \times \widehat{S}_d$ that satisfy $\underline{V} \leq \overline{V}$ by construction, and $\overline{V} \leq \underline{V}$ by the comparison principle. By the uniqueness result for viscosity solutions to the HJB equation \cite[Theorem 9]{mfcp}, this implies that $ \overline{V} = \underline{V} = \widehat{V}$, and the construction of $\overline{V}$ and $\underline{V}$ provides uniform convergence as a result. 
\end{itemize}

\subsection{Notation}
As noted above, we define $\lbrs{d} := \{1, \ldots, d\}$. The $(d-1)$--dimensional simplex $S_d$ is given by
\begin{align}
    \label{eq:simplex}
    S_d := \Big\{(m_1, \ldots, m_d) \in \R^{d} : m_i \geq 0, \quad \sum_{i=1}^d m_i = 1\Big\}.
\end{align}
On the simplex $S_d$, we define directional derivatives for a function $\varphi : S_d \to \R$ by
\begin{align*}
    D^i f(m) := (\partial_{m_j - m_i} f(m))_{j \in \lbrs{d}},
\end{align*}
where derivatives are in the direction of $e_j - e_i$ for the corresponding standard basis vectors of $\R^d$.

Unless otherwise specified, $\Omega$ will refer to a subset of $\R^n$ (measurable, but not necessarily open). As usual, $L^p(\Omega)$ consists of the measurable functions on $\Omega$ such that
\begin{align*}
    \|f\|_{p} := \left(\int_{\Omega} |f(x)|^p dx\right)^{1/p} < \infty,
\end{align*}
where the integral above is taken with respect to the Lebesgue measure on $\R^d$, restricted to $\Omega$. For $p = \infty$, we instead take $L^\infty(\Omega)$ to be the space of measurable functions on $\Omega$ whose essential supremum is bounded. We say that a function is differentiable on $\Omega$ if it is differentiable on the interior of $\Omega$ and extends to a differentiable function on an open neighborhood of $\R^n$ containing $\Omega$. As is standard, $\sC(\Omega)$ refers to the space of continuous functions on $\Omega$, equipped with the uniform norm:
\begin{align*}
    \|f\|_{\infty} := \sup_{x \in \Omega} |f(x)|.
\end{align*}
Below, when we refer to \textit{uniform convergence}, we always refer to convergence in the above norm.

We use the notation $\sC^{1,1}(\Omega)$ to denote (continuously-differentiable) functions that have Lipschitz-continuous first derivatives on $\Omega$. See \cite{mfcp} for a more thorough discussion of this class of functions as they relate to the HJB equation in Equation~\eqref{eq:mfcp-hjb}. Other spaces of functions that we refer to throughout the paper include the space of continuously-differentiable functions with $k$ continuously-differentiable derivatives, denoted by $\sC^k(\Omega)$. We define a norm on $\sC^k(\Omega)$ by
\begin{align*}
    \|f\|_{\sC^k(\Omega)} := \max_{|\alpha| \leq m} \sup_{x \in \Omega} |D^\alpha f(x)|.
\end{align*}

\subsection{Organization}
We structure the remainder of the paper as follows. Section~\ref{sec:hjb-mfcp-setup} introduces the relevant assumptions placed upon the MFCP that ensure sufficient regularity of the value function. Section~\ref{sec:dgm-conv} presents our modified version of the DGM algorithm for the HJB equation in Equation~\eqref{eq:mfcp-hjb}, the corresponding convergence results, and all associated proofs. Finally, Section~\ref{sec:numerical-tests} presents selected numerical results from a JAX implementation of the DGM.
\section{Assumptions and Known Results for the MFCP}
\label{sec:hjb-mfcp-setup}
Below, we outline the (already known) MFCP and the recent existence, uniqueness, and convergence results for the MFCP proven in \cite{mfcp}. We also present the relevant assumptions imposed upon the MFCP. This section utilizes the background on stochastic control and HJB equations provided in~\cite{stoch-control-text}.

\subsection{Assumptions for the MFCP}
\label{subsec:assumptions}
We start by assuming the following, drawing upon some of the notations introduced throughout the preceding section.
\begin{assumption*}{Assumption A}
\label{a:mfcp-A}
The function $F : [0, T] \times [0, M]^{d^2} \times S_d \to \R$, given by
\begin{align*}
    F(t, a^1, \ldots, a^d, m) := \sum_{i \in \lbrs{d}} m_i f(t, i, a^i, m)
\end{align*}
is continuous and there exists a constant $C>0$, such that for any $t,s\in[0,T]$, $a\in[0,M]^d$, and $m,p\in S^d$,
\begin{align*}
    |F(t,a,m) - F(s,a,p)| &\leq  C(|t-s| + |m-p|).
\end{align*}
\end{assumption*}
\begin{assumption*}{Assumption B}
\label{a:mfcp-B}
For each $i \in \lbrs{d}$, the running cost $f$ is continuously differentiable in $a$, $\nabla_a f$ is Lipschitz continuous with respect to $m$, and there exists $\lambda > 0$ such that
\begin{align*}
    f(t, i, b,m) \geq f(t, i, a,m) + \langle \nabla_a f(t, i, a, m), b-a \rangle + \lambda |b-a|^2.
\end{align*}
\end{assumption*}
\begin{assumption*}{Assumption C}
\label{a:mfcp-C}
We have that:
\begin{enumerate}
\item[(C1)] Defining $G : S_d \to \R$ by
\begin{align*}
    G(m) := \sum_{i \in \lbrs{d}} m_i g^i(m),
\end{align*}
we have that
$F(\cdot, a, \cdot) \in \sC^{1,1}([0, T] \times S_d)$ and $G \in \sC^{1,1}(S_d)$.
\item[(C2)] The function
\begin{align*}
    [0,T] \times [0, \infty)^{d \times d} \times \intr(S_d) \ni (t,w,m) \mapsto \sum_{i \in \lbrs{d}}m_i f\Big(t, i, \Big(\frac{w_{i,j}}{m_i}\Big)_{j\neq i}, m\Big)
\end{align*}
is convex in $(w, m)$.
\item[(C3)] $G$ is convex in $m$.
\end{enumerate}
\end{assumption*}
Note that under Assumption~\hyperref[a:mfcp-C]{(C)}, $G$ is Lipschitz continuous: for any $m, p \in S_d$,
\begin{align*}
    |G(m) - G(p)| &\leq C|m-p|
\end{align*}
for some $C > 0$, which we can take to be the same constant as in Assumption~\hyperref[a:mfcp-A]{(A)}.

We provide an example of an MFCP in Section~\ref{sec:numerical-tests} for which the running and terminal costs satisfy all three of the above assumptions.
\subsection{Known Convergence, Existence, and Uniqueness results for the MFCP}
\label{subsec:prev-results}
Under assumptions similar to those in the preceding section, Cecchin \cite{mfcp} derived a series of useful results for the MFCP. To begin, Equation~\eqref{eq:mfcp-hjb} has a unique classical solution $V \in \sC^{1,1}([0,T] \times S_d)$.
\begin{proposition}[{{\cite[Theorem 2.9(3)]{mfcp}}}]
\label{prop:mfcp-hjb}
Under Assumptions~\hyperref[a:mfcp-A]{(A)}--~\hyperref[a:mfcp-C]{(C)} presented in Section~\ref{subsec:assumptions} above, the value function $V \in \sC^{1,1}([0, T] \times S_d)$ for the MFCP is the unique classical solution of the HJB equation \eqref{eq:mfcp-hjb}, with Hamiltonian defined in Equation~\eqref{eq:hamiltonian}.
\end{proposition}

\begin{remark}
\normalfont
If we only assume that Assumptions~\hyperref[a:mfcp-A]{(A)} and ~\hyperref[a:mfcp-B]{(B)} hold, then $V$ is both Lipschitz continuous in $(t, m)$ and the unique viscosity solution of Equation~\eqref{eq:mfcp-hjb} on $S_d$, and an optimal control for the MFCP exists \cite[Theorem 2.9(2)]{mfcp}. Under Assumption~\hyperref[a:mfcp-A]{(A)}, Cecchin also presents several convergence results connecting the MFCP to the $N$-agent optimization problem. In particular, we have, from \cite[Theorem 2.10]{mfcp}, the convergence of the value function $V^N$ for the $N$-agent optimization problem in \eqref{eq:n-agent} to the value function $V$ for the MFCP with a rate of $O(1/\sqrt{N})$.  Furthermore, \cite[Theorem 2.11]{mfcp} establishes the existence of an asymptotic optimal control for the $N$-agent optimization problem using the optimal control for the MFCP, also with an $O(1/\sqrt{N})$ rate of convergence.
\end{remark}

If in addition to Assumption~\hyperref[a:mfcp-A]{(A)}, we also include Assumptions~\hyperref[a:mfcp-B]{(B)}, then the empirical distribution corresponding to the optimal control for the MFCP, given in Equation~\eqref{eq:optimal-control-mfcp}, converges to the optimal process for the $N$-agent problem, with a convergence rate of $O(N^{-1/9})$ \cite[Theorem 2.13]{mfcp}.

\section{The DGM Algorithm and the Main Convergence Results}
\label{sec:dgm-conv}
In this section, we present our modified DGM algorithm to solve Equation \eqref{eq:mfcp-hjb} and the corresponding convergence proof. We begin by presenting our version of the DGM algorithm in Section~\ref{sec:dgm}. Next, in Section~\ref{sec:nn-approximators}, we describe the relevant results from universal approximation theory in addition to proving that the DGM can approximate the solution to Equation~\eqref{eq:mfcp-hjb} arbitrarily well in Theorem~\ref{th:unif-error-conv}. In Section~\ref{sec:nn-convergence}, we present our main convergence results in Theorem~\ref{th:unif-value-conv}. Finally, in Section~\ref{subsec:comp-original-dgm}, we 
describe the difficulties in establishing convergence for the algorithm with $L^2$-loss and explain with more details why the proof of \cite{dgm} is invalid for our equation \eqref{eq:mfcp-hjb}. We employ the following outline to obtain our existence and convergence results:
\begin{enumerate}
\item[(1)] By an appropriate version of the universal approximation theorem, we can arbitrarily approximate $V \in \sC^{1,1}([0,T] \times S_d)$ with neural network approximators; see Proposition~\ref{th:univ-approx}. Specifically, there exists a sequence of neural network approximators $\{\varphi(t, m; \theta^n)\}_{n \in \N}$ such that the DGM loss applied to the sequence converges to zero as $n \to \infty$; see Theorem~\ref{th:unif-error-conv}. 
\item[(2)] If the DGM loss goes to zero along this sequence of neural network approximators $\{\varphi(t, m; \theta^n)\}_{n \in \N}$, then $\varphi(\cdot,\cdot;\theta^n) \to V(\cdot,\cdot)$ uniformly on $[0, T] \times S_d$; see Theorem~\ref{th:unif-value-conv}.
\end{enumerate}
Throughout this section, we utilize the following operator, defined for $\phi \in \sC^{1,1}([0,T] \times S_d)$:
\begin{align}\label{eq:operator1}
    \sL[\phi](t,m) := -\partial_t \phi(t, m) + \sum_{i \in \lbrs{d}} m_i H^i(t, m, D^i \phi(t, m)).
\end{align}
Note that the first line in Equation \eqref{eq:mfcp-hjb} can be written as $\sL[V](t,m)=0$.
\subsection{DGM Algorithm}
\label{sec:dgm}
In this section, we present a modification of the DGM algorithm, first proposed by \cite{dgm}, in the context of the HJB equation for the MFCP. 
The DGM aims to efficiently approximate a solution to the above equation using a deep learning-based approach. Specifically, the method learns model parameters $\theta \in \R^P$, where $P$ depends on the dimension $d$ of the simplex $S_d$, the number of layers in the neural network in use, and the number of nodes in each layer of the neural network. The DGM learns the model parameters $\theta$ by minimizing an objective functional, referred to as the {\it DGM loss}, and given by 
\begin{align}
\label{eq:unif-loss}
\begin{split}
L(\theta) &:= \max_{(t, m) \in [0, T] \times S_d} |\sL[\varphi(\cdot, \cdot; \theta)](t,m)| + \max_{m \in S_d}|\varphi(T, m;\theta) - \sum_{i \in \lbrs{d}}m_i g^i(m)|. 
\end{split}
\end{align}
On occasions, we refer to this loss as the $L^\infty$-loss. To this end, we refer to our algorithm as the {\it uniform DGM}, referring to the fact that we use an approximation of the $L^\infty$-loss functional defined above to obtain uniform convergence to the true value function.
The maxima over $[0, T] \times S_d$ and $S_d$ are approximated by sampling, as demonstrated in the following algorithm. As with the original DGM algorithm with $L^2$-loss \cite{dgm}, we utilize stochastic gradient descent (SGD) to find the parameter $\theta \in \R^P$ that minimizes an appropriately defined loss. Note that in Algorithm~\ref{alg:unif-dgm}, the architecture of the neural network is fixed, and only the parameter $\theta$ is updated by SGD.

Importantly, the $L^\infty$-loss in Equation~\eqref{eq:unif-loss} is intractable to compute exactly. Instead, given a set of $M$ samples $(t^{(j)}, m^{(j)}, p^{(j)})_{j = 1, \ldots, M}:=\{(t^{(j)}, m^{(j)}, p^{(j)}):j = 1, \ldots, M\}  \subset [0,T] \times S_d\times S_d$, we define
\begin{align}
\label{eq:sampled-unif-loss}
\begin{split}
    G((t^{(j)}, m^{(j)}, p^{(j)})_{j = 1, \ldots, M}, \theta) &:= \max_{j = 1, \ldots, M}|\sL[\varphi(\cdot, \cdot; \theta)](t^{(j)},m^{(j)})| \\
    & \ + \max_{j = 1, \ldots, M} \Big|\varphi(T, p^{(j)};\theta) - \sum_{i \in \lbrs{d}}p_i^{(j)} g^i(p^{(j)})\Big|.
\end{split}
\end{align}
If $\rho$ is the uniform distribution on $[0, T] \times S_d \times S_d$, then our algorithm attempts to minimize the following population loss via SGD:
\begin{align}
\label{eq:pop-unif-loss}
\begin{split}
    L_M(\theta) &:= \E_{(t^{(j)}, m^{(j)}, p^{(j)})_{j = 1, \ldots, M} \sim \rho^{\otimes M}} \bigg[\max_{j = 1, \ldots, M}|\sL[\varphi(\cdot, \cdot; \theta)](t^{(j)},m^{(j)})| \\
    &\qquad + \max_{j = 1, \ldots, M} \Big|\varphi(T, p^{(j)};\theta) - \sum_{i \in \lbrs{d}}p_i^{(j)} g^i(p^{(j)})\Big|\bigg] \\
    &= \E_{(t^{(j)}, m^{(j)}, p^{(j)})_{j = 1, \ldots, M} \sim \rho^{\otimes M}} [G((t^{(j)}, m^{(j)}, p^{(j)})_{j = 1, \ldots, M}, \theta)].
\end{split}
\end{align}
Note, in particular, that if $(t^{(j)}, m^{(j)}, p^{(j)})_{j = 1, \ldots, M}$ are drawn uniformly, then $G((t^{(j)}, m^{(j)}, p^{(j)})_{j = 1, \ldots, M}, \theta)$ is an unbiased estimate of $L_M(\theta)$.
\begin{algorithm}[H]
\caption{Uniform DGM}\label{alg:unif-dgm}
\begin{algorithmic}[1]
\State Initialize parameters $\theta^{(0)} \in \R^P$
\State Initialize tolerance $\delta \in (0,1)$
\State $n \gets 0$
\State Sample $(t^{(j)}, m^{(j)}, p^{(j)})_{j = 1, \ldots, M} \in [0,T] \times S_d$
\While{$G((t^{(j)}, m^{(j)}, p^{(j)})_{j = 1, \ldots, M}, \theta^{(n)}) \geq \delta$} \Comment{See Equation~\eqref{eq:sampled-unif-loss}}
\State Sample $(t^{(j)}, m^{(j)}, p^{(j)})_{j = 1, \ldots, M} \in [0,T] \times S_d$
\State $\theta^{(n+1)} \gets \theta^{(n)} - \alpha^{(n)} \nabla_\theta G((t^{(j)}, m^{(j)}, p^{(j)})_{j = 1, \ldots, M}, \theta^{(n)})$
\State $n \gets n + 1$
\EndWhile
\end{algorithmic}
\end{algorithm}
Moreover, the stochastic gradient $\nabla_\theta G((t^{(j)}, m^{(j)}, p^{(j)})_{j = 1, \ldots, M}, \theta)$ is also an unbiased estimate of $L_M(\theta)$, as is required in order to access the convergence guarantees associated with SGD \cite[Chapter 6]{bubeck-co}.

As shown in \cite[Lemma 1]{cohen2024deep}, we have that
\begin{align*}
    \E|L(\theta) - G((t^{(j)}, m^{(j)}, p^{(j)})_{j = 1, \ldots, M}, \theta)| \leq C_d M^{-1/(d-1)},
\end{align*}
where the constant $C_d > 0$ depends on the dimension $d$ and the parameters of Equation~\eqref{eq:mfcp-hjb}, and the expectation is taken with respect to the samples $(t^{(j)}, m^{(j)}, p^{(j)})_{j = 1, \ldots, M}$, drawn uniformly from $[0, T] \times S_d \times S_d$. In turn, we have that
\begin{align*}
    |L(\theta) - L_M(\theta)| &= \left|L(\theta) - \E[G((t^{(j)}, m^{(j)}, p^{(j)})_{j = 1, \ldots, M}, \theta)]\right| \\
    &\leq \E|L(\theta) - G((t^{(j)}, m^{(j)}, p^{(j)})_{j = 1, \ldots, M}, \theta)| \\
    &\leq C_d M^{-1/(d-1)}.
\end{align*}
If we have found parameters $\theta \in \R^{P}$ such that $L_M(\theta) < \varepsilon / 2$ for some $\varepsilon > 0$, and $M$ is taken sufficiently large so that $C_d M^{-1/(d-1)} < \varepsilon / 2$, we can conclude that $L(\theta) < \varepsilon$, a condition that is sufficient to access our convergence guarantee in Theorem~\ref{th:unif-value-conv}. More explicitly, this discussion can be summarized in the following proposition, justifying the construction of Algorithm~\ref{alg:unif-dgm}.

\begin{proposition}
\label{prop:algorithmic-guarantee}
Suppose $\varepsilon > 0$, and assume that Algorithm~\ref{alg:unif-dgm} returns parameters $\theta^* \in \R^P$ such that the population loss in Equation~\eqref{eq:pop-unif-loss} satisfies $L_M(\theta^*) < \varepsilon / 2$. If $M > (2 C_d / \varepsilon)^{d - 1}$, then the $L^\infty$-loss in Equation~\eqref{eq:unif-loss} satisfies $L(\theta^*) < \varepsilon$.
\end{proposition}

In practice, the performance of Algorithm~\ref{alg:unif-dgm} may vary depending on the sample size $M$ at each step. Additionally, the learning rate schedule $\alpha^{(n)}$ may determine the convergence rate of the algorithm as before. Using an optimizer such as AdaGrad or Adam may help speed up convergence. Finally, instead of using a tolerance $\delta \in (0,1)$ to determine the convergence of the algorithm, one may instead specify a fixed number of SGD iterations to carry out.

The structure of the PDE in Equation~\eqref{eq:mfcp-hjb} prohibits us from using the same argument as \cite{dgm} to prove convergence of the DGM algorithm with $L^2$-loss. Instead, by slightly changing the loss function to the $L^\infty$-loss, given in \eqref{eq:unif-loss}, which we use to train the neural network approximation in the DGM algorithm, we can prove the convergence of our modified DGM algorithm to the unique value function of the MFCP. Specifically, by utilizing the above loss functional that approximates the uniform norm of the PDE and terminal condition rather than the squared error of the PDE and the terminal condition, we can leverage the theory of viscosity solutions for first-order HJB equations from \cite{visc-users} in our convergence proof. 
\begin{remark}
\normalfont
Note that we do not discuss convergence of the optimization method used in Algorithm 1, and we do not establish that it is possible to find the global minimizer of the loss functional using Algorithm~\ref{alg:unif-dgm}. Rather, by ``establishing convergence of the uniform DGM,'' we intend to show that if one can minimize the loss functional using an appropriate training method, then the resulting neural network approximators will converge uniformly to the unique solution to Equation~\eqref{eq:mfcp-hjb}. Of course, in practice, a method such as Algorithm~\ref{alg:unif-dgm} is by far the most common method for training neural networks to achieve this goal. Additionally, as seen in Proposition~\ref{prop:algorithmic-guarantee}, if the population loss defined in Equation~\eqref{eq:pop-unif-loss} is made small via SGD, we can ensure that the $L^\infty$-loss in Equation~\eqref{eq:unif-loss} is also small, assuming enough samples are taken at each step of training.
\end{remark}
Finally, it is worth noting that the maximum is generally not a smooth function. In practice, the maximum may be computed using a smooth maximum -- a function that approximates a maximum but is differentiable. For previous work involving stochastic gradient descent applied to a maximized loss function and insights into its robustness, refer to \cite{pmlr-v48-shalev-shwartzb16}. 
\subsection{Neural Network Approximators}
\label{sec:nn-approximators} 
This section starts by utilizing the classical universal approximation result of Hornik \cite{univ-approx-thm} to establish that the uniform DGM loss can be made arbitrarily small by neural network approximators in Theorem~\ref{th:unif-error-conv}. Then, we lay the foundations for our main convergence result via Corollary~\ref{corr:convergence-to-v}, which provides a reformulation of the result in Theorem~\ref{th:unif-error-conv} that is better suited for the language of viscosity solutions. 

Although all universal approximation results in \cite{univ-approx-thm} hold for two-layer neural networks (i.e., a single hidden layer), we allow neural networks with multiple layers, including deep neural networks with modern architectures. In general, a network with $L$ layers, maximum width $n$, and a common activation function $\sigma$ takes the form
\begin{align}
\label{eq:multi-layer-nn}
    \varphi(t, m; \theta) := \sigma(W_{L} \ldots \sigma(W_1 m + \alpha t + c_1) \ldots + c_{L}),
\end{align}
where the activation function $\sigma$ is applied elementwise. Above, $W_i$ are weight matrices, $c_i$ are bias vectors, and $\alpha$ is a scalar weight. In turn, the parameters of each neural network are of the form $\theta = (W_1, \ldots, W_L, c_1, \ldots, c_L, \alpha) \in \R^P$ (upon flattening all weight matrices into vectors), where $P$ depends on the maximum width of the network, the depth $L$ of the network, and the dimension $d$ of the simplex $S_d$. In turn, we take $\mathfrak{C}^{(P)}_{d+1}(\sigma)$ to be the class of neural networks with parameters $\theta$ of dimension at most $P$, from which we define
\begin{align*}
    \mathfrak{C}_{d+1}(\sigma) := \bigcup_{P=1}^\infty \mathfrak{C}^{(P)}_{d+1}(\sigma).
\end{align*}
This is a slight departure from the notation of \cite{univ-approx-thm}, but all the relevant universal approximation results therein still hold in this more general context.
With the above notation in mind, we can apply the universal approximation theorem \cite[Theorem 3]{univ-approx-thm}, stated as follows:
\begin{proposition}
\label{th:univ-approx}
If $\sigma \in \sC^m(\R)$ is nonconstant and bounded, then $\mathfrak{C}_{d+1}$ is uniformly $m$-dense on compact sets in $\sC^m(\R^{d+1})$. In particular, for all $h \in \sC^m(\R^{d+1})$, all compact subsets $K \subset \R^{d+1}$, and any $\varepsilon > 0$, there exists $\psi = \psi(h, K, \varepsilon) \in \mathfrak{C}_{d+1}$ such that $\|h - \psi\|_{\sC^m(K)} < \varepsilon$.
\end{proposition}
In the implementation in Section~\ref{sec:numerical-tests}, we take $\sigma(y) = \tanh(y)$, a typical choice of activation function that is smooth, nonconstant, and bounded -- and therefore satisfies all the criteria of Proposition~\ref{th:univ-approx}. Note also that with this choice of $\sigma$, any element of $\mathfrak{C}_{d+1}(\sigma)$ is smooth (as a linear combination of smooth functions), ensuring that any $\varphi \in \mathfrak{C}_{d+1}(\sigma)$ has Lipschitz-continuous first derivative. See \cite{tanh-approx} for further justification of this choice of activation function in terms of the approximation guarantees that it brings.

With the above background in mind, we move towards approximating solutions to the HJB equation for the MFCP \eqref{eq:mfcp-hjb} using the DGM. Through assumptions~~\hyperref[a:mfcp-A]{(A)} --~\hyperref[a:mfcp-C]{(C)} in \cite{mfcp}, we have the following useful properties for proving convergence of the DGM to the solution of the HJB equation.
\begin{enumerate}
\item[(1)] Under Assumption~\hyperref[a:mfcp-A]{(A)}, the PDE-Hamiltonian in Equation~\eqref{eq:mfcp-hjb}, namely, $H(t, m, z)=\sum_{i \in \lbrs{d}} m_i H^i(t, m, z)$, is Lipschitz continuous in all three arguments. Specifically, we have that
\begin{align*}
    \Big|\sum_{i \in \lbrs{d}} m_i(H^i(t_1, m_1, z_1) - H^i(t_2, m_2, z_2))\Big| &\leq \kappa|(t_1 - t_2, m_1 - m_2, z_1 - z_2)| \\
    &\leq \kappa(|t_1 - t_2| + |m_1 - m_2| + |z_1 - z_2|)
\end{align*}
for some constant $\kappa > 0$, where $|\cdot|$ denotes the Euclidean norm. In fact, we only need the following, which also follows from the assumptions presented in \cite{mfcp}: for each $i \in \lbrs{d}$ and fixed $(t,m) \in [0,T] \times S_d$, the map $p \mapsto H^i(t, m, p)$ is Lipschitz continuous with a common Lipschitz constant $C > 0$.
\item[(2)] By Proposition \ref{prop:mfcp-hjb}, Equation~\eqref{eq:mfcp-hjb} admits a unique classical solution $V \in \sC^{1,1}([0,T] \times S_d)$.
\end{enumerate}

We show below that by utilizing the uniform error, given in Equation~\eqref{eq:unif-loss} 
we can obtain the desired convergence result. With the ultimate goal of showing that a neural network can approximate the value function $V(t, m)$ on $[0, T] \times S_d$ arbitrarily well in the uniform norm by taking %$n$, 
the number of neurons in the network, sufficiently large, we first show the following {\it existence result}. This result in fact holds for both the our DGM loss functional and the DGM with $L^2$-loss. In particular, the following theorem establishes the existence of a sequence of neural networks that makes the DGM loss arbitrarily small.
\begin{theorem}
\label{th:unif-error-conv}
Let $\sigma \in \sC^1(\R)$ be bounded and nonconstant. For every $\varepsilon > 0$, there exists a constant $\Tilde{K}(d, T, C) > 0$, where $d$ is the dimension of the simplex $S_d$, $T$ is the finite time horizon of the MFCP, and $C$ is the Lipschitz constant of the PDE-Hamiltonian in Equation~\eqref{eq:mfcp-hjb}, such that for some $\varphi=\varphi(\cdot,\cdot;\theta) \in \mathfrak{C}_{d+1}(\sigma)$, the DGM loss functional in Equation~\eqref{eq:unif-loss} satisfies $L(\theta) \leq \Tilde{K} \varepsilon$.
\end{theorem}
\begin{proof}
Note that $\Omega_T := [0, T] \times S_d$ is a compact set in $\R^{d+1}$. Thus, by Proposition~\ref{th:univ-approx} above, we know that for the unique solution $V \in \sC^{1,1}(\Omega_T)$ to Equation~\eqref{eq:mfcp-hjb} and any $\varepsilon > 0$, there exists $\varphi \in \mathfrak{C}_{d+1}(\sigma)$ such that
\begin{align}
    \label{eq:univ-approx-est}
    \begin{split}
    \varepsilon > \sup_{(t,m) \in \Omega_T} |V(t,m) - \varphi(t, m; \theta)| &+ \sup_{(t,m) \in \Omega_T} |\partial_t V(t,m) - \partial_t \varphi(t, m; \theta)| \\
    &+ \sup_{(t, m) \in \Omega_T} |\nabla_m V(t,m) - \nabla_m \varphi(t,m; \theta)|  
    \end{split}
\end{align}
For such $\varphi \in \mathfrak{C}_{d+1}(\sigma)$, the Lipschitz continuity of the PDE-Hamiltonian in Equation~\eqref{eq:mfcp-hjb} yields
\begin{align*}
    \Big|\sum_{i \in \lbrs{d}} m_i H^i(t, m, D^i \varphi(t, m; \theta)) - \sum_{i \in \lbrs{d}} m_i H^i(t, m, D^i V(t, m))\Big| \leq C \sum_{i \in \lbrs{d}}|D^i \varphi(t, m; \theta) - D^iV(t,m)|,
\end{align*}
where $C$ is the Lipschitz constant of the PDE-Hamiltonian. Above, we use the fact that $|m_i| \leq 1$ for any $m \in S_d$. Next, for each $i \in \lbrs{d}$ and any $(t, m) \in \Omega_T$, we observe that, by denoting the standard basis of $\R^d$ by $\{e_i\}_{i \in \lbrs{d}}$, we have that
\begin{align}
\label{eq:gradient-error-estimate}
    |D^i \varphi(t, m; \theta) - D^iV(t,m)|^2 &= \frac{1}{2}\sum_{j=1}^d |(\nabla_m \varphi(t,m; \theta)- \nabla_m V(t,m)) \cdot (e_j - e_i)|^2 \notag \\
    &\leq \frac{1}{2}\sum_{j=1}^d |\nabla_m V(t,m) - \nabla _m\varphi(t,m; \theta)|^2 |e_j - e_i|^2 \notag \\
    &\leq d |\nabla_m V(t,m) - \nabla_m \varphi(t,m; \theta)|^2,
\end{align}
applying the Cauchy--Schwarz inequality in the second line above. Putting both bounds together, we obtain
\begin{align}
\label{eq:unif-error-estimate}
    \Big|\sum_{i \in \lbrs{d}} m_i H^i(t, m, D^i \varphi(t, m; \theta)) - \sum_{i \in \lbrs{d}} m_i H^i(t, m, D^i V(t, m))\Big| \leq C d^{3/2} |\nabla_m V(t,m) - \nabla_m \varphi(t,m; \theta)|.
\end{align}
Now, let $\Tilde{K} = C d^{3/2} > 0$. We conclude by noting that the value function $V$ satisfies $\sL[V](t,m) = 0$ for all $(t,m) \in \Omega_T$ (in addition to the terminal condition of Equation~\eqref{eq:mfcp-hjb}), which allows us to conclude that for $\varphi=\varphi(\cdot,\cdot;\theta)$,
\begin{align*}
    L(\theta) &= \max_{(t,m) \in \Omega_T}\left|\sL[\varphi](t,m)\right| + \max_{m \in S_d}\Big|\varphi(T, m; \theta) - \sum_{i \in \lbrs{d}} m_i g^i(m)\Big| \\
    &= \max_{(t,m) \in \Omega_T}\left|\sL[\varphi](t,m)\ - \sL[V](t,m)\right|+ \max_{m \in S_d}\left|\varphi(T, m; \theta) - V(T, m)\right|\\
    &\leq \max_{(t,m) \in \Omega_T}\Big|m_i H^i(t, m, D^i \varphi(t, m; \theta)) - \sum_{i \in \lbrs{d}} m_i H^i(t, m, D^i V(t, m))\Big| \\
    &\quad + \max_{(t,m) \in \Omega_T}\left|\partial_t V(t,m) - \partial_t \varphi(t,m;\theta)\right|\\
    &\quad + \max_{m \in S_d}\left|\varphi(T, m; \theta) - V(T, m)\right|\\
    &\leq \Tilde{K} \varepsilon
\end{align*}
by applying the approximation result from \eqref{eq:univ-approx-est} and taking $\Tilde{K} > 0$ larger if necessary. Note that the constant $\tilde{K}$ now depends on $d, C,$ and $T$ as claimed.
\end{proof}

We include the constant $\tilde{K}(d, T, C) > 0$ in front of $\varepsilon$ in the statement of Theorem 3.3 to emphasize that the loss functional's value depends not only on the approximation error $\varepsilon >0 $ that arises from Proposition~\ref{th:univ-approx} but also on $d$, $T$, and the Lipschitz constant of the Hamiltonian. This constant draws a distinction between simply approximating the value function directly using a neural network and attempting to minimize the loss functional with a neural network approximator.
\begin{remark}
\normalfont
We note here that the proof of Theorem~\ref{th:unif-error-conv} is quite similar to the proof of \cite[Theorem 7.1]{dgm}, the analogous result from the original paper on the DGM. Our general proof technique is essentially the same, although the individual computations differ substantially because we use $L^\infty$-loss instead of $L^2$-loss, and the original DGM paper's theoretical analysis considers a class of quasilinear parabolic PDEs in divergence form.
\end{remark}
The following corollary, which also relates to the existence of an approximating sequence of neural networks, utilizes the same universal approximation theorem as Theorem~\ref{th:unif-error-conv}. In particular, we can obtain a sequence of neural networks that satisfies a corresponding sequence of PDEs with a measurable error term that uniformly converges to zero.
\begin{corollary}
\label{corr:convergence-to-v}
There exists a sequence of 
parameters $\{\theta^n\}_{n \in \N}$, such that $\varphi^n \to V$ uniformly as $n \to \infty$, where $V$ is the unique classical solution to Equation~\eqref{eq:mfcp-hjb},
and $\varphi^n(t,m):=\varphi(t,m;\theta^n)$.
Given such $\varphi^n$, we define $e^n : [0, T] \times S^d \to \R$ as follows:
\begin{align}
\label{eq:corr-hjb}
\begin{split}
    &e^n(t, m):=\sL[\varphi^n](t, m) ,\qquad t\in[0,T),\\
    &e^n(T, m):=\varphi^n(T, m) - \sum_{i \in \lbrs{d}} m_i g^i(m) .
\end{split}
\end{align}
Then, $\|e^n\|_{\infty} \to 0$ as $n \to \infty$, noting that each $e^n$ depends on the neural network parameter $\theta^n$.
\end{corollary}
\begin{proof}
From Proposition~\ref{th:univ-approx}, we know that there exists a sequence of neural networks $\{\varphi^n\}_{n \in \N} \subset \mathfrak{C}_{d+1}(\sigma) \subset \sC^{1,1}([0,T] \times S_d)$, parametrized by a set of parameters $\{\theta^n\}_{n \in \N}$, such that \eqref{eq:univ-approx-est} holds for $\varepsilon = n^{-1}$. Thus, it immediately follows that $\|\varphi^n - V\|_{\infty} \to 0$ as $n \to \infty$,  yielding a sequence $\varphi^n$ of neural networks parametrized by $\theta^n$ that converges uniformly to $V$, the classical solution to Equation~\eqref{eq:mfcp-hjb}.
By construction, $\varphi^n$ satisfies the PDE in \eqref{eq:corr-hjb}. Furthermore, because each $e^n: [0,T] \times S_d \to \R$ is continuous on both $[0, T) \times S_d$ and $\{T\} \times S_d$, it is evidently measurable.

Finally, to see that $\|e^n\|_\infty \to 0$ as $n \to \infty$, we reuse many of the estimates from the proof of Theorem~\ref{th:unif-error-conv}. In particular, we can write
\begin{align*}
    e^n(t, m) &= -\partial_t \varphi^n(t, m) + \sum_{i \in \lbrs{d}} m_i H^i(t, m, D^i \varphi^n(t,m)) \\
    &= \partial_t V(t,m) -\partial_t \varphi^n(t, m) + \sum_{i \in \lbrs{d}} m_i\left[H^i(t, m, D^i \varphi^n(t,m))- H^i(t, m, D^i V(t,m))\right],
\end{align*}
for $(t,m) \in [0,T) \times S_d$, using the fact that $V$ solves Equation~\eqref{eq:mfcp-hjb}. By the fact that $m \in S_d$ and the Lipschitz continuity of $H^i$ for each $i \in \lbrs{d}$, we again have that
\begin{align*}
    \Big|\sum_{i \in \lbrs{d}} m_i\left[H^i(t, m, D^i \varphi^n(t,m))- H^i(t, m, D^i V(t,m))\right]\Big| &\leq C \sum_{i \in \lbrs{d}} |D^i \varphi^n(t, m) - D^iV(t,m)| \\
    &\leq 2 dC |\nabla_m V(t,m) - \nabla_m \varphi^n(t,m)|.
\end{align*}
Thus, we have that for all $(t,m) \in [0,T) \times S_d$
\begin{align*}
    |e^n(t,m)| \leq 2dC|\nabla_m V(t,m) - \nabla_m \varphi^n(t,m)| + |\partial_t V(t,m) -\partial_t \varphi^n(t, m)| \leq \frac{2dC}{n} + \frac{1}{n},
\end{align*}
by the construction of $\varphi^n$ from \eqref{eq:univ-approx-est}. Similarly, at the terminal time $t = T$, observe that
\begin{align*}
    |e^n(T,m)| = \Big|\varphi^n(T, m) - \sum_{i \in \lbrs{d}} m_i g^i(m)\Big| = |\varphi^n(T,m) - V(T, m)| \leq \frac{1}{n}.
\end{align*}
Thus, it follows that $\|e^n\|_{\infty} \to 0$ as $n \to \infty$ as claimed.
\end{proof}
\begin{remark}
\normalfont
The above sequence of neural networks satisfies several properties based on its construction. Namely, the sequence $\{\varphi^n\}_{n \in \N}$ is uniformly bounded and equicontinuous on $[0, T] \times S_d$. This is a consequence of the converse of the Arzelà--Ascoli theorem. For details, see \cite[Theorem 7.24]{rudin-pma}, which shows that if a sequence of continuous functions on a compact metric space converges uniformly, then the sequence must be equicontinuous. Moreover, because our sequence is uniformly convergent, it is sequentially compact and hence bounded, given that we are working in a metric space. This property is a consequence of the construction of our loss functional, although it does not play a further role in our results.
\end{remark}
\subsection{Convergence of Neural Network Approximators to Value Function}
\label{sec:nn-convergence}
We now discuss the convergence of a sequence of neural network approximators $\varphi^n$ to $V$, the unique classical solution of the HJB equation in Equation \eqref{eq:mfcp-hjb}. With the ultimate goal of establishing \textit{uniform} convergence of the neural network approximators $\varphi^n$ to the value function $V$, we rely on the theory of viscosity solutions to first-order nonlinear PDEs. This powerful theory, developed by Crandall, Evans, and Lions in the 1980s for the explicit purpose of approaching HJB equations (which often lack classical, differentiable solutions) \cite{visc-users}, will allow us to relate the neural network approximators $\varphi^n(t,m)$ to the value function $V$ via a sequence of first-order nonlinear PDEs. Then, using a version of the comparison principle for viscosity solutions, we obtain the desired convergence.

Below, it will be useful to instead consider the HJB equation for $\eta \in \widehat{S}_d$, where
\begin{align*}
    \widehat{S}_d := \Big\{(\eta_1, \ldots, \eta_{d-1}) \in \R^{d-1} : \eta_j \geq 0 \; \text{for all} \; j = 1,\ldots, d-1, \; \sum_{j=1}^{d-1} \eta_j \leq 1 \Big\}.
\end{align*}
Then, the simplex $S_d$ can be expressed as
\begin{align*}
    S_d = \Big\{(\eta, \eta^{-d}) : \eta \in \widehat{S}_d, \; \eta^{-d} = 1 - \sum_{j=1}^{d-1} \eta_j\Big\}.
\end{align*}
In turn, any function $v \in \sC^1(S_d)$ induces a function $\widehat{v} \in \sC^1(\widehat{S}_d)$, given by $  \widehat{v}(\eta)=v(\eta, \eta^{-d})$. Furthermore, the gradient $\nabla_\eta \widehat{v}(\eta) = (\partial_{\eta_j}\widehat{v}(\eta))_{j=1}^{d-1}$ satisfies $\partial_{\eta_j}\widehat{v}(\eta) = \partial_{m_j - m_d} v(m)$. Following \eqref{eq:operator1}, we define an operator for functions defined on  $[0,T]\times \widehat{S}_d$ by
\begin{align}
\label{eq:mfcp-operator}
    \widehat{\sL}[\widehat{\phi}](t,\eta) := -\partial_t \widehat{\phi}(t, \eta) + \sum_{i \in \lbrs{d-1}} \eta_i \widehat{H}^i(t, \eta, \nabla_\eta \widehat{\phi}(t,\eta)) + \eta^{-d} \widehat{H}^d(t, \eta, \nabla_\eta \widehat{\phi}(t,\eta)).
\end{align}
Now, the solution to Equation~\eqref{eq:mfcp-hjb} can be written $V(t, \eta, \eta^{-d}) = \widehat{V}(t, \eta)$, where $\widehat{V} \in \sC^{1,1}(\widehat{S}_d)$ is the unique solution to the modified HJB equation 
\begin{align}
\label{eq:chart-hjb}
\begin{split}
    &\widehat{\sL}[\widehat{V}](t,\eta)= 0,\\
    &\widehat{V}(T, \eta) = \sum_{i \in \lbrs{d-1}} \eta_i g^i(\eta, \eta^{-d}) + \eta^{-d}g^d(\eta, \eta^{-d}).
\end{split}
\end{align}
Above, the modified Hamiltonians take inputs in $[0, T] \times \widehat{S}_d \times \R^{d-1}$ and are defined by
\begin{align*}
    \widehat{H}^i(t, \eta, p) = H^i(t, \eta, \eta^{-d}, p_1 - p_i, \ldots, p_{d-1} - p_i, -p_{i}) \quad \text{ and } \quad \widehat{H}^d(t, \eta, p) = H^d(t, \eta, \eta^{-d}, p, 0),
\end{align*}
for $i\in\lbrs{d-1}$. 
As shown in \cite{mfcp}, Equation~\eqref{eq:chart-hjb} has a unique solution $\widehat{V} \in \sC^{1,1}([0,T] \times \widehat{S}_d)$. Additionally, $\widehat{S}_d$ is a compact subset of $\R^{d-1}$, allowing us to apply the universal approximation theorem exactly as above, now on $\widehat{S}_d$. By the definitions of the modified Hamiltonians (in terms of the original Hamiltonians), we can also apply the exact same argument as above to obtain a result equivalent to our Theorem~\ref{th:unif-error-conv} on $\widehat{S}_d$. In particular, a straightforward computation shows that for all $i \in \lbrs{d}$ and $p,p' \in \widehat{S}_d$,
\begin{align*}
    |\widehat{H}^i(t, \eta, p) - \widehat{H}^i(t, \eta, p')| \leq D |p - p'|,
\end{align*}
where $D = \sqrt{2(2(d-1) + 1)} C$. This shows that $\widehat{H}^i$ is Lipschitz continuous in $p$ with Lipschitz constant $D \geq C > 0$ for $i \in \lbrs{d}$. From this, the proof of a modified version of Theorem~\ref{th:unif-error-conv}, now on $[0, T] \times \widehat{S}_d$, can proceed exactly as before. Note that the original value function defined on the simplex can be recovered via $V(t, \eta, \eta^{-d}) = \widehat{V}(t, \eta)$ for $\eta \in \widehat{S}_d$. 

Now, we can reframe the problem in terms of Equation~\eqref{eq:chart-hjb}, which possesses a unique classical solution $\widehat{V} \in \sC^{1,1}([0, T] \times \widehat{S}_d)$. Recall that $\widehat{S}_d \subset \R^{d-1}$ is the preimage of the simplex in $\R^d$ under the chart introduced in Equation~\eqref{eq:chart-hjb}. Working with Equation~\eqref{eq:chart-hjb} rather than Equation~\eqref{eq:mfcp-hjb} allows us to cite results from the theory of viscosity solutions that require the domain of the relevant PDE to be open; note that $\intr(\widehat{S}_d)$ is an open subset of $\R^{d-1}$ whereas $S_d$ has empty interior in $\R^d$. Additionally, the following result demonstrates that the convergence result on $\widehat{S}_d$ translates to $S_d$ without any issues.
\begin{proposition}
\label{prop:eps-n-convergence}
Assume that $\widehat{V} \in \sC([0,T] \times \widehat{S}_d)$ and $\widehat{\phi}^n \in \sC([0,T] \times \widehat{S}_d)$ are such that $\|\widehat{V} - \widehat{\phi}^n\|_\infty \to 0$ as $n \to \infty$. Then, $\|V - \phi^n\|_\infty \to 0$ as $n \to \infty$, where $V, \phi^n \in \sC([0,T] \times S_d)$ are given by $V(t, \eta, \eta^{-d}) = \widehat{V}(t, \eta)$ and $\phi^n(t, \eta, \eta^{-d}) = \widehat{\phi}^n(t, \eta)$ for all $(t, \eta, \eta^{-d}) \in [0,T] \times S_d$ and all $n \in \N$.
\end{proposition}
We remark that the opposite direction is also true; given a sequence of functions on $[0, T] \times S_d$ that converge uniformly to $V \in \mathcal{C}([0, T] \times S_d)$, the corresponding functions on $[0, T] \times \widehat{S}_d$ converge to $\widehat{V} \in \mathcal{C}([0, T] \times \widehat{S}_d)$, given by $\widehat{V}(t, \eta) = V(t, \eta, \eta^{-d})$ in our notation above. However, we do not require the converse of Proposition~\ref{prop:eps-n-convergence}, and the proof is analogous.
\begin{proof}
This is a simple consequence of the definition of $\widehat{S}_d$. Indeed, if $\|\widehat{V} - \widehat{\phi}^n\|_\infty \to 0$ as $n \to \infty$. Then, for any $\varepsilon > 0$, we have that for all $n \in \N$ sufficiently large,
\begin{align*}
    \sup_{(t,\eta) \in [0,T] \times \widehat{S}_d} |\widehat{V}(t, \eta) - \widehat{\phi}^n(t, \eta)| < \varepsilon.
\end{align*}
Consequently, for all $n \in \N$ sufficiently large, we have that
\begin{align*}
    \sup_{(t, \eta, \eta^{-d}) \in [0,T] \times S_d} |V(t,\eta, \eta^{-d}) - \phi^n(t,\eta, \eta^{-d})| &= \sup_{(t,\eta) \in [0,T] \times \widehat{S}_d} |\widehat{V}(t,\eta) - \widehat{\phi}^n(t,\eta)| \\
    &< \varepsilon.
\end{align*}
This means precisely that $\|V - \phi^n\|_\infty \to 0$ as $n \to \infty$.
\end{proof}
As a consequence of the above proposition, it suffices to show the uniform convergence of a sequence of neural network approximators $\widehat{\varphi}^n$ to the unique classical solution $\widehat{V}$ of Equation~\eqref{eq:chart-hjb} on $\widehat{S}_d$, as we can then recover uniform convergence on the simplex. Note that from now on, we consider neural networks, denoted by $\widehat{\varphi}$, which are defined as in \eqref{eq:multi-layer-nn} but with the input $\eta\in\widehat{S}_d$ rather than $m\in S_d$. From the discussion preceding Equation~\eqref{eq:chart-hjb}, we also know that Theorem~\ref{th:unif-error-conv} holds on $\widehat{S}_d$, yielding the existence of a sequence of neural networks $\{\widehat{\varphi}^n(t,\eta):=\widehat{\varphi}(t,\eta; \theta^n)\}_{n \in \N}$ such that
$L(\theta^n) \to 0$ as $n \to \infty$. In turn, each network $\widehat{\varphi}^n(t,\eta)$ satisfies its own ``perturbed'' PDE, of the form
\begin{align}
\label{eq:perturbed-hjb}
    \widehat{\sL}[\widehat{\varphi}^n](t, \eta) = \widehat{e^n}(t, \eta), \qquad (t,\eta) \in [0,T] \times \widehat{S}_d,
\end{align}
with $\widehat{e}^n(t, \eta):=e^n(t,\eta,\eta^{-d})$. 
For notational simplicity, we take
\begin{align*}
    \widehat{G}(\eta) := \sum_{i \in \lbrs{d-1}} \eta_i g^i(\eta, \eta^{-d}) + \eta^{-d}g^d(\eta, \eta^{-d})
\end{align*}
in this section to denote the terminal condition of the HJB equation on $\widehat{S}_d$. Denoting $\widehat{\Omega}_T := [0, T] \times \widehat{S}_d$ as in the previous section, Theorem~\ref{th:unif-error-conv} above implies that Equation~\eqref{eq:perturbed-hjb} satisfies
\begin{align*}
    \max_{(t,\eta) \in [0,T] \times \widehat{S}_d}\left|\widehat{e^n}(t,\eta)\right|+ \max_{\eta \in \widehat{S}_d}\left|\widehat{\varphi}^n(T, \eta) - \widehat{G}(\eta)\right| \to 0\qquad\text{as}\qquad n\to\infty.
\end{align*}
With this context in mind, we state the main {\it convergence result} of this section.
\begin{theorem}
\label{th:unif-value-conv}
The family of neural network approximators $\{\widehat{\varphi}^n(t, \eta)\}_{n \in \N}$ satisfying Equation~\eqref{eq:perturbed-hjb} converges uniformly to $\widehat{V} \in \sC^{1,1}([0,T] \times \widehat{S}_d)$, the unique classical solution of Equation~\eqref{eq:chart-hjb}, in the sense that
\begin{align*}
    \sup_{(t,\eta) \in [0,T] \times \widehat{S}_d} |\widehat{\varphi}^n(t, \eta) - \widehat{V}(t, \eta)| \to 0\qquad \text{as}\qquad n\to\infty.
\end{align*}
\end{theorem}
To prove the above theorem, we argue via the comparison principle for viscosity solutions to \eqref{eq:mfcp-hjb} presented in \cite{mfcp}. To this end, we require a suitable definition of viscosity solutions of Equation~\eqref{eq:chart-hjb} on $\widehat{S}_d$.
\begin{definition}
\label{def:visc-solutions}
A function $\widehat{v} \in \sC((0, T) \times \intr(\widehat{S}_d))$ is: 
\begin{enumerate}
\item[(i)] a viscosity subsolution of Equation~\eqref{eq:chart-hjb} if for any $\widehat{\varphi} \in \sC^1((0, T) \times \intr(\widehat{S}_d))$, $\widehat{\mathcal{L}}[\widehat{\varphi}] \leq 0$
for every local maximum $(t_0, \eta_0) \in (0, T) \times \intr(\widehat{S}_d)$ of $\widehat{v} - \widehat{\varphi}$ on $(0, T) \times \intr(\widehat{S}_d)$.
\item[(ii)] a viscosity supersolution of Equation~\eqref{eq:chart-hjb} if for any $\widehat{\varphi} \in \sC^1((0, T) \times \intr(\widehat{S}_d))$, $\widehat{\mathcal{L}}[\widehat{\varphi}] \geq 0$
for every local minimum $(t_0, \eta_0) \in (0, T) \times \intr(\widehat{S}_d)$ of $\widehat{v} - \widehat{\varphi}$ on $(0, T) \times \intr(\widehat{S}_d)$.
\item[(iii)] a viscosity solution of Equation~\eqref{eq:chart-hjb} if $\widehat{v}$ is both a viscosity subsolution and viscosity supersolution.
\end{enumerate}
\end{definition}
\begin{remark}
\normalfont
When viscosity solutions are introduced in \cite{mfcp}, the author also allows for test functions on $[0, T) \times S_d$ (resp. $[0, T) \times \widehat{S}_d$), noting that $[0, T) \times S_d$ (resp. $[0, T) \times \widehat{S}_d$) is no longer an open subdomain of $\R^{d+1}$ (resp. $\R^d$). However, in order to utilize \cite[Theorem 3.3]{visc-users}, the standard comparison principle for viscosity solutions, we must consider viscosity solutions on open subdomains of $\R^{d}$. As noted in \cite{mfcp}, it is also not immediately clear that a classical solution to Equation~\eqref{eq:mfcp-hjb} is a viscosity solution if the latter is defined on a closed set. 

We could alternatively cite the comparison principle from \cite[Theorem 3.4]{mfcp} that utilizes the definition of viscosity solutions on closed sets presented therein. However, in order to utilize the clearly presented stability properties of viscosity solutions under uniform limits presented in \cite{ crandall-evans-lions, visc-users}, we opt for the standard definition in Definition~\ref{def:visc-solutions}.
\end{remark}
In order to establish Theorem~\ref{th:unif-value-conv}, we proceed using a standard comparison principle argument for viscosity solutions that also leverages the fact that $\widehat{V} \in \sC^{1,1}([0, T] \times \widehat{S}_d)$ is the unique viscosity solution to Equation~\eqref{eq:chart-hjb} from \cite{mfcp}.
\begin{proof}[Proof of Theorem~\ref{th:unif-value-conv}]
For each $n \in \N$, we may define an operator
\begin{align*}
    \widehat{\sL}^n[\widehat{\phi}](t, \eta) := -\partial_t \widehat{\phi}(t, \eta) + \sum_{i \in \lbrs{d-1}} \eta_i \widehat{H}^i(t, \eta, \nabla_\eta \widehat{\phi}(t,\eta)) + \eta^{-d} \widehat{H}^d(t, \eta, \nabla_\eta \widehat{\phi}(t,\eta)) - \widehat{e^n}(t,\eta),
\end{align*}
corresponding to the sequence of PDEs described in \eqref{eq:perturbed-hjb}. Because $\widehat{\sL}^n[\widehat{\phi}]$ depends only on the derivatives of $\widehat{\phi}$ (and not on $\widehat{\phi}$ itself), we observe that $\widehat{\sL}^n$ is \textit{proper} in the sense of \cite{visc-users}.\footnote{In \cite{visc-users}, the authors results only apply to PDEs of the form $F(x, u, \nabla u, \nabla^2 u) = 0$, where $F : \R^n \times \R \times \R^n \times \R^{n \times n} \to \R$ satisfies the monotonicity condition $F(x, r, p, X) \leq F(x, s, p, Y)$ when $r \leq s$ and $X \preceq Y$ in the Loewner order on symmetric matrices. If $F$ satisfies condition, as is the case for the class of HJB equations that we consider, then the corresponding PDE is proper. Our setting belongs under the umbrella of degenerate elliptic equations, discussed in \cite{visc-users}, that do not depend on the second derivatives of $u$.} This fact also ensures that the technical conditions preceding the comparison principle \cite[Theorem 3.3]{visc-users} are satisfied.

Now, note that by the discussion following Proposition~\ref{prop:eps-n-convergence},
$$\max_{(t,\eta) \in [0,T] \times \widehat{S}_d}\left|\widehat{e^n}(t,\eta)\right| \to 0\qquad \text{as}\qquad n\to\infty,$$
meaning that $\widehat{e^n}$ converges uniformly to zero on $[0,T] \times \widehat{S}_d$. Now, for each $n \in \N$, define $T^n : [0,T] \times \widehat{S}_d \times \R \times \R^{d-1} \to \R$ by
\begin{align*}
    T^n(t, \eta, p_0, p) := -p_0 + \sum_{i \in \lbrs{d-1}} \eta_i \widehat{H}^i(t, \eta, p) + p^{-d} \widehat{H}^d (t,\eta,p) - \widehat{e^n}(t, \eta).
\end{align*}
We then have that $T^n$ converges uniformly on $[0,T] \times \widehat{S}_d \times \R \times \R^{d-1}$ to
\begin{align*}
    T(t,\eta,p_0, p) := -p_0 + \sum_{i \in \lbrs{d-1}} \eta_i \widehat{H}^i(t, \eta, p) + p^{-d} \widehat{H}^d (t,\eta,p).
\end{align*}
These definitions are motivated by the fact that Equation~\eqref{eq:mfcp-hjb} can be written succinctly as
\begin{align*}
    T(t,\eta, \partial_t \widehat{\varphi}, \nabla_\eta \widehat{\varphi}) = 0,
\end{align*}
while Equation~\eqref{eq:perturbed-hjb} is given by
\begin{align*}
    T^n(t,\eta, \partial_t \widehat{\varphi}, \nabla_\eta \widehat{\varphi}) = 0
\end{align*}
for each $n \in \N$. Now, following \cite[Remark 6.3]{visc-users}, we note that because $\widehat{\varphi}^n$ is a classical solution (and therefore a viscosity solution) to the equation $T^n(t,\eta, \partial_t \widehat{\varphi}, \nabla_\eta \widehat{\varphi}) = 0$ on $(0,T) \times \intr(\widehat{S}_d)$, then 
\begin{align*}
    \overline{V}(t, \eta) := \lim_{j \to \infty} \sup\{\widehat{\varphi}^n(s, \nu) : n \geq j, \; (t, \eta) \in [0,T) \times \widehat{S}_d, \; |(s, \nu) - (t, \eta)| \leq 1/j\}
\end{align*}
is a viscosity subsolution to the equation $T(t,\eta, \partial_t \widehat{\varphi}, \nabla_\eta \widehat{\varphi}) = 0$, as we have that
\begin{align*}
    T(t,\eta,p_0,p) = \liminf_{n \to \infty} T^n(t,\eta,p,p_0).
\end{align*}
On the other hand, we also observe that 
\begin{align*}
    \underline{V}(t, \eta) := \lim_{j \to \infty} \inf\{\widehat{\varphi}^n(s,\nu) : n \geq j, \; (t, \eta) \in [0,T) \times \widehat{S}_d, \; |(s, \nu) - (t,\eta)| \leq 1/j\}
\end{align*}
is a viscosity supersolution to the equation $T(t,\eta, \partial_t \widehat{\varphi}, \nabla_\eta \widehat{\varphi}) = 0$ by the same reasoning. By construction, observe that $\underline{V} \leq \overline{V}$ on $[0,T) \times \widehat{S}_d$. Note also that both $\underline{V}$ and $\overline{V}$ are well-defined on $\{0\} \times \partial \widehat{S}_d$ by their construction. However, by the comparison principle presented in \cite[Theorem 3.3]{visc-users}, the fact that $\underline{V}$ is a viscosity supersolution and $\overline{V}$ is a viscosity subsolution is sufficient to conclude that $\overline{V} \leq \underline{V}$ on $[0,T) \times \widehat{S}_d$, observing that the comparison principle still holds on the closure of the domain $(0, T) \times \intr(\widehat{S}_d)$.

In particular, $\overline{V} = \underline{V}$ is a viscosity solution. As shown in \cite[Theorem 9]{mfcp}, Equation~\eqref{eq:mfcp-hjb} has a unique viscosity solution $\widehat{V}$, showing that $\overline{V} = \underline{V} = \widehat{V}$. Now, \cite[Remark 6.4]{visc-users} implies that $\lim_{n \to \infty} \widehat{\varphi}^n(t, \eta) = \widehat{V}(t, \eta)$ uniformly on $[0, T) \times \widehat{S}_d$. Finally, note that we also have that
\begin{align*}
   \max_{\eta \in \widehat{S}_d}\left|\widehat{\varphi}^n(T, \eta) - \widehat{V}(T, \eta)\right| \to 0,
\end{align*}
from the construction of the modified DGM loss, allowing us to conclude uniform convergence of $\widehat{\varphi}^n \to \widehat{V}$ on the entire region $[0, T] \times \widehat{S}_d$ as claimed.
\end{proof}

\subsection{Comparison to the DGM with $L^2$-loss}
\label{subsec:comp-original-dgm}
\normalfont
At this point, we can clarify the reasons for the modification to the DGM algorithm made in Section~\ref{sec:dgm}. Here, we refer to the Sobolev space $W^{k,p}(\Omega)$, given by the space of functions $f \in L^p(\Omega)$ such that for any multi-index $\alpha$ with $|\alpha| \leq k$, the derivative $D^\alpha f$ exists and belongs to $L^p(\Omega)$ itself; see \cite[Chapter 5]{evans-pde} for more details on this characterization of Sobolev spaces. We also define the space $L^p(0, T; W^{k, q}(\Omega))$ by the function $f \in L^p([0, T] \times \Omega)$ such that, for fixed $t \in [0, T]$, $f(t, \cdot) \in W^{k, q}(\Omega)$. The spaces $L^p(0, T; L^q(\Omega))$ are then defined analogously for $1 \leq p, q \leq \infty$. 

Sirignano and Spiliopoulos \cite{dgm} formulated the $L^2$-loss,\footnote{As with the DGM loss in Equation~\eqref{eq:unif-loss}, this loss is computed in practice by sampling points according to probability measures $\nu_1$ and $\nu_2$ on $[0, T] \times S_d$ and $S_d$ respectively to obtain an unbiased estimate of the $L^2$-error: in Equation~\eqref{eq:l2-error}} as
\begin{align}
\label{eq:l2-error}
    \Tilde{L}(\theta) = \|\sL[\varphi(\cdot, \cdot; \theta)](t,m)\|_{2, [0, T] \times S_d, \nu_1}^2 + \|\varphi(T, m; \theta) - \sum_{i \in \lbrs{d}} m_i g^i(m)\|_{2, S_d, \nu_2}^2,
\end{align}
because of the natural connection between the class of equations that they considered and convergence in $L^2$. A key step in the proof of their analogue to Theorem~\ref{th:unif-value-conv} involves obtaining a uniform bound on $\{\varphi(t, m; \theta^n)\}_{n \in \N}$ in $L^\infty(0, T; L^2(\Omega)) \cap L^2(0, T; W^{1,2}(\Omega))$, where $\Omega$ is the open domain on which the PDE is considered. In turn, this arises from an energy bound on quasilinear parabolic equations such as the one presented in \cite{energy-est-orig}, or in more generality in \cite{parabolic-est}. However, such a bound only holds for equations of the form
\begin{align*}
\begin{cases}
    \partial_t u - \text{div}(a(t,x, u, \nabla u)) = H(t,x, \nabla u) & (t,x) \in (0,T) \times \Omega, \\
    u = 0 & (t,x) \in (0,T) \times \partial \Omega, \\
    u(0,x) = u_0(x) & x \in \Omega.
\end{cases}
\end{align*}
that satisfy the Leray--Lions conditions. Namely, there must exist $\alpha > 0$ such that
\begin{align*}
    \alpha|\xi|^p \leq a(t,x, p, \xi) \cdot \xi
\end{align*}
for all $\xi \in \R^d$ and some $1 < p < d$. Clearly, this fails in our case, where $a$ is identically zero, even though our HJB equation otherwise satisfies the structure conditions in \cite{parabolic-est}.

Translating the convergence argument in \cite[Theorem 7.3]{dgm} to our context is also complicated by the fact that the class of quasilinear parabolic PDEs for which they prove convergence of the DGM possesses a standard notion of weak solutions that, via the dominated convergence theorem, cooperates with convergence in $L^2$. In the case of HJB equations, however, viscosity solutions take the place of weak solutions and instead behave nicely with respect to uniform convergence. Thus, we require that $\widehat{e^n} \to 0$ uniformly on $[0, T] \times \widehat{S}_d$ in Equation~\eqref{eq:perturbed-hjb}, but the formulation of the DGM with $L^2$-loss only implies convergence of the error term in $L^2$.

Finally, \cite{dgm} only concludes uniform convergence of the neural network approximators to the true solution of the PDE after imposing additional assumptions of uniform boundedness and equicontinuity on the neural networks. Although this may be a reasonable assumption to include, we find that our modified DGM algorithm and the theory of viscosity solutions provide a more direct route to uniform convergence.

Here, we also note an interesting connection between our work and recent work in the field of physics-informed neural networks (PINNs), a recent framework for incorporating PDE constraints into the training of neural networks. Specifically \cite{l2-pinns} show that, for a class of second-order HJB equations arising from stochastic control problems, it is \textit{impossible} to obtain a convergence guarantee along the lines of Theorem~\ref{th:unif-value-conv} if $L^2$-loss is used. We emphasize that this result does not necessarily apply to the first-order HJB equation arising from the MFCP. Our work is further complicated by the fact that Equation~\eqref{eq:mfcp-hjb} has solutions on the simplex $[0, T] \times S_d$, whereas the class of stochastic control problems considered in \cite{l2-pinns} are solved over $[0, T] \times \R^n$. However, future work may investigate whether the stability results obtained in \cite{l2-pinns} generalize to the context of the MFCP. Currently, this connection only provides some intuition as to why the $L^\infty$-loss may be the appropriate choice for training neural networks to solve Equation~\eqref{eq:mfcp-hjb}.
\section{Numerical Results}\label{sec:numerical-tests}
In this section, we present numerical results for the DGM applied to a simple example case of the MFCP, as presented in \cite[Example 2]{mfcp}. In particular, we work with the following example.
\begin{example}
\normalfont
Consider the quadratic running cost
\begin{align*}
    f(t, i, \alpha, m) = \frac{1}{2} \sum_{j\in\lbrs{d}, j \neq i} c_{i,j} \alpha_{i,j}^2 + f_0^i(m),
\end{align*}
with $f_0^i(m):= m_i$ and $\{c_{i,j}\}_{i,j \in \lbrs{d}} \in \R^{d \times d}$ a cost matrix that encodes the cost of transitioning from state $i$ to state $j$ for $i \neq j$. Finally, we consider the linear terminal condition given by $g^i(m) = m_i$ for $i \in \lbrs{d}$. With this choice of terminal cost, we obtain the terminal condition
\begin{align*}
    V(T, m) = G(m) = \sum_{i \in \lbrs{d}} m_i^2.
\end{align*}
In this simple example, the Hamiltonian is explicitly given 
\begin{align*}
    H^i(t,m,z) = \sum_{j \neq i} \left(-\mathfrak{a}^*\left(-\frac{z_j}{c_{i,j}}\right)z_j - \frac{1}{2}\left(\mathfrak{a}^*\left(-\frac{z_j}{c_{i,j}}\right)\right)^2\right) - f_0^i(m),
\end{align*}
where
\begin{align*}
    \mathfrak{a}^*(s) = \begin{cases}
    0 & s \leq 0, \\
    s & 0 \leq s \leq M, \\
    M & s \geq M.
    \end{cases}
\end{align*}
Recall that $M > 0$ is some constant such that $A = [0, M]^d$, where $A$ is the action space for the MFCP. Under this construction, all of the convexity and Lipschitz continuity constraints in Assumptions~\hyperref[a:mfcp-A]{(A)} -- \hyperref[a:mfcp-C]{(C)} are satisfied.

In the context of this example setup, we compare the performance of the DGM with $L^2$-loss from \cite{dgm} and the modified algorithm with our $L^\infty$-loss. Additionally, we explore the dependence of the performance of Algorithm~\ref{alg:unif-dgm} on the number of samples drawn from $[0, T] \times S_d \times S_d$. Finally, we demonstrate that the modified DGM algorithm scales well with dimension, providing numerical tests up to dimension $d = 200$. All numerical experiments were run on a single NVIDIA A100 Tensor Core GPU via Google Colaboratory, and all code can be found in the \href{https://github.com/jakehofgard/jax-dgm}{GitHub repository} for this paper. 

For the sake of comparison, we utilize the same architecture as in \cite{dgm}; a neural network with three hidden LSTM-like layers and one dense output layer, with $\tanh$ activation throughout. We performed hyperparameter tuning on the number of layers, the number of samples used in each epoch, the width of the network, and the learning rate schedule. Unless otherwise noted, all results were produced by training the neural network for 10000 epochs, with 10 gradient steps in each epoch, and $M = 5000$ samples for each epoch. After hyperparameter tuning, we found that a learning rate cosine decay schedule, with a linear warmup period, achieved the best performance. All experiments were performed with the Adam optimizer with weight decay, which we found improved performance. Finally, we found that multiplying the second term in the loss functional in Equation~\eqref{eq:sampled-unif-loss} by a positive constant $\lambda_n > 0$ (where the subscript $n$ indicates that this constant may change as Algorithm~\ref{alg:unif-dgm} progresses) often had a positive effect on the optimization procedure. Note that the presence of such a constant in the loss functional has no effect on our convergence proof in the preceding section. Adaptively decreasing this constant from a large value such $\lambda_n = 100$ to $\lambda_n = 2$ within 5000 epochs encourages the algorithm to learn the terminal condition more quickly, although this minor change to Algorithm~\ref{alg:unif-dgm} is not strictly necessary.

\subsection{Results and Discussion}
In dimension $d = 2$, it is possible to compute a precise approximation of the true solution to the HJB equation using a finite difference scheme. We implemented such a solver, using an explicit Euler scheme on the unit interval. Indeed, using the formulation in Equation~\eqref{eq:chart-hjb}, the HJB equation in $d = 2$ can be viewed as a backwards parabolic PDE on the unit interval. Our implementation utilizes the \texttt{py-pde} package \cite{Zwicker2020}, a package designed for solving nonlinear parabolic PDEs on spatial grids such as the unit interval. In Figure~\ref{fig:d=2-true-sol}, we display the true solution obtained using our finite difference solver.

On the other hand, the plots in Figure~\ref{fig:d=2-unifsurface} and Figure~\ref{fig:d=2-l2surface} below demonstrate the value functions approximated by DGM with $L^\infty$-loss and $L^2$-loss respectively in dimension $d=2$. Additionally, the two plots display the absolute error between the learned solutions and the finite difference solution in Figure~\ref{fig:d=2-true-sol}. Figures~\ref{fig:l2-loss} and~\ref{fig:unif-loss} contain the corresponding loss curves. We remark that the loss curves for the two methods are displayed on different scales because they are \textit{not} comparable metrics of the performance of the solver. Given the formulation of the $L^\infty$-loss in Equation~\ref{eq:unif-loss}, we expect the $L^\infty$-loss to converge to a higher value than the $L^2$-loss, even if both metrics result in models that approximately solve the HJB equation. We note that even though the $L^\infty$-loss stabilizes at a higher value than the $L^2$-loss, the resulting model appears much more stable \textit{and} solves the HJB equation to a higher degree of accuracy according to the error between learned solution and the finite difference solution. Specifically, the model trained with $L^\infty$-loss obtains a maximum absolute error of $0.055$, a mean absolute error of $0.022$, and a mean squared error of $0.00067$. Meanwhile, the model trained with $L^2$-loss obtains a maximum absolute error of approximately $0.11$, a mean absolute error of $0.024$, and mean square error of $0.0011$. This disparity reinforces the notion that training with $L^\infty$-loss encourages better convergence.

\begin{figure}[htbp]
\centering
  \includegraphics[width=0.5\linewidth]{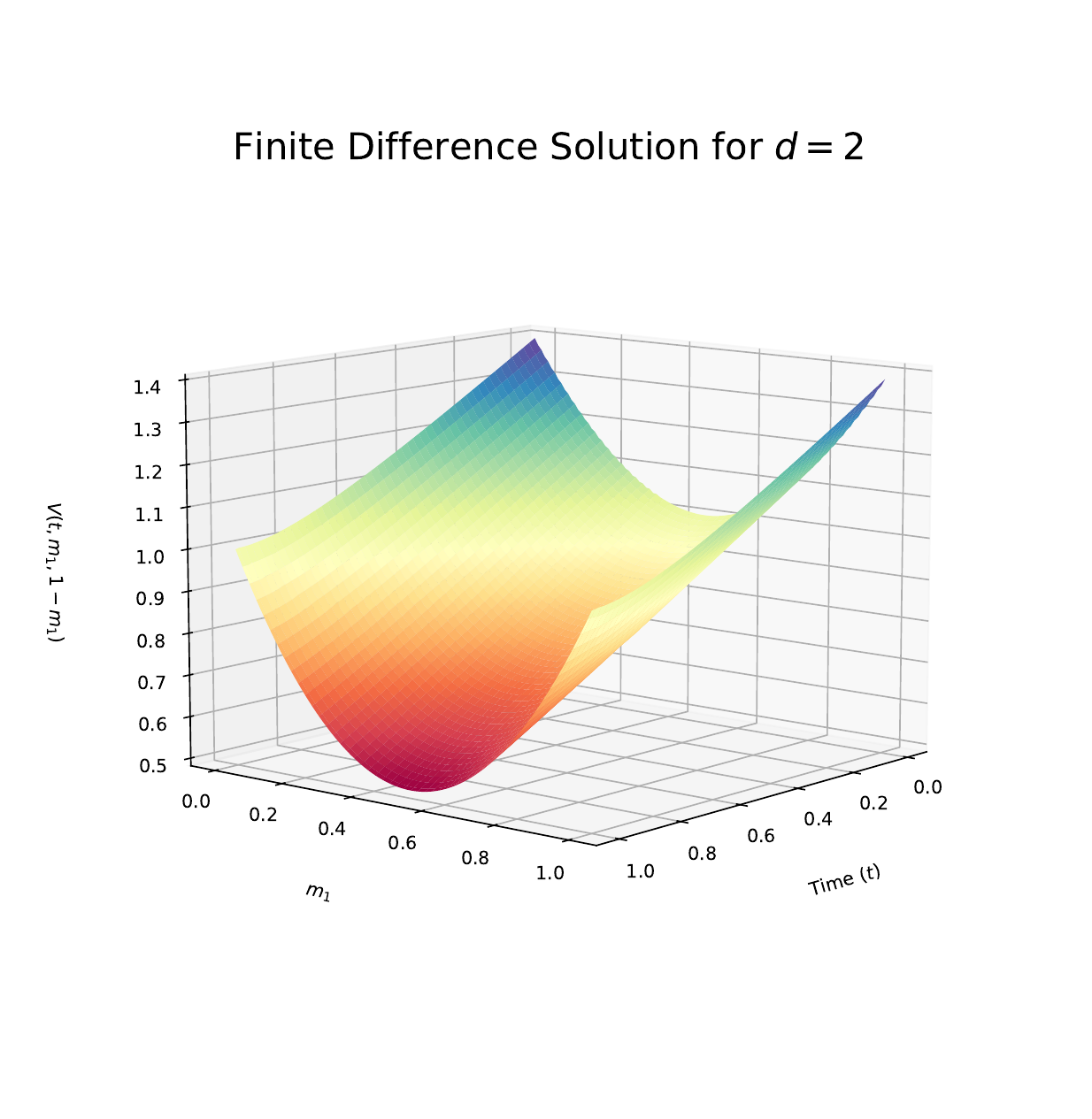}
  \captionof{figure}{\footnotesize True solution for $d = 2$, computed using a explicit finite difference scheme.}
  \label{fig:d=2-true-sol}
\end{figure}

\begin{figure}[htbp]
\centering
\begin{minipage}{0.49\textwidth}
  \centering
  \includegraphics[width=\linewidth]{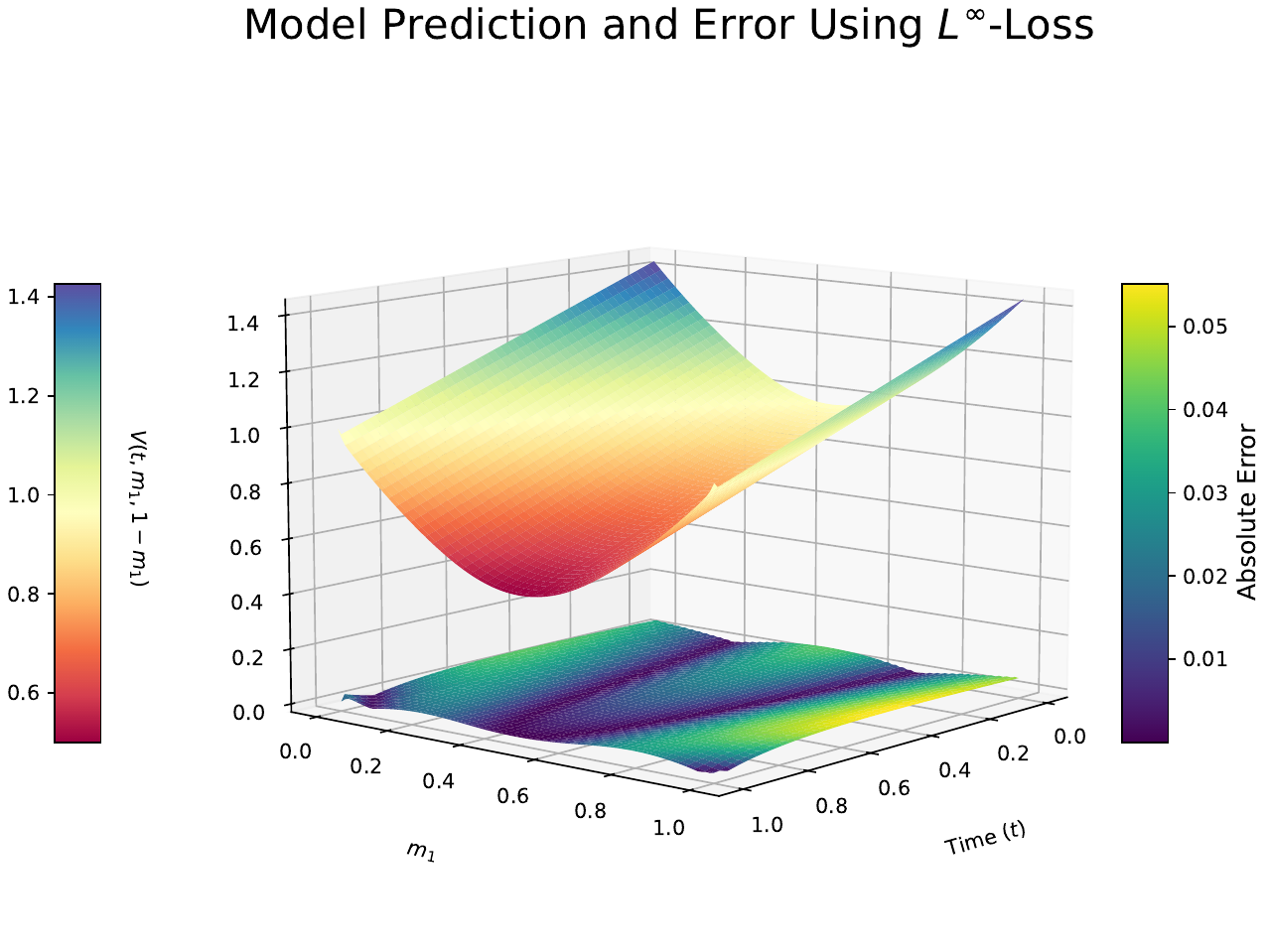}
  \captionof{figure}{\footnotesize Approximate value function, learned using $L^\infty$-loss, in dimension $d = 2$. The absolute error, measured between the approximate value function and the finite difference solution in Figure~\ref{fig:d=2-true-sol} is displayed on the same plot.}
  \label{fig:d=2-unifsurface}
\end{minipage}
\begin{minipage}{0.49\textwidth}
  \centering
  \includegraphics[width=\linewidth]{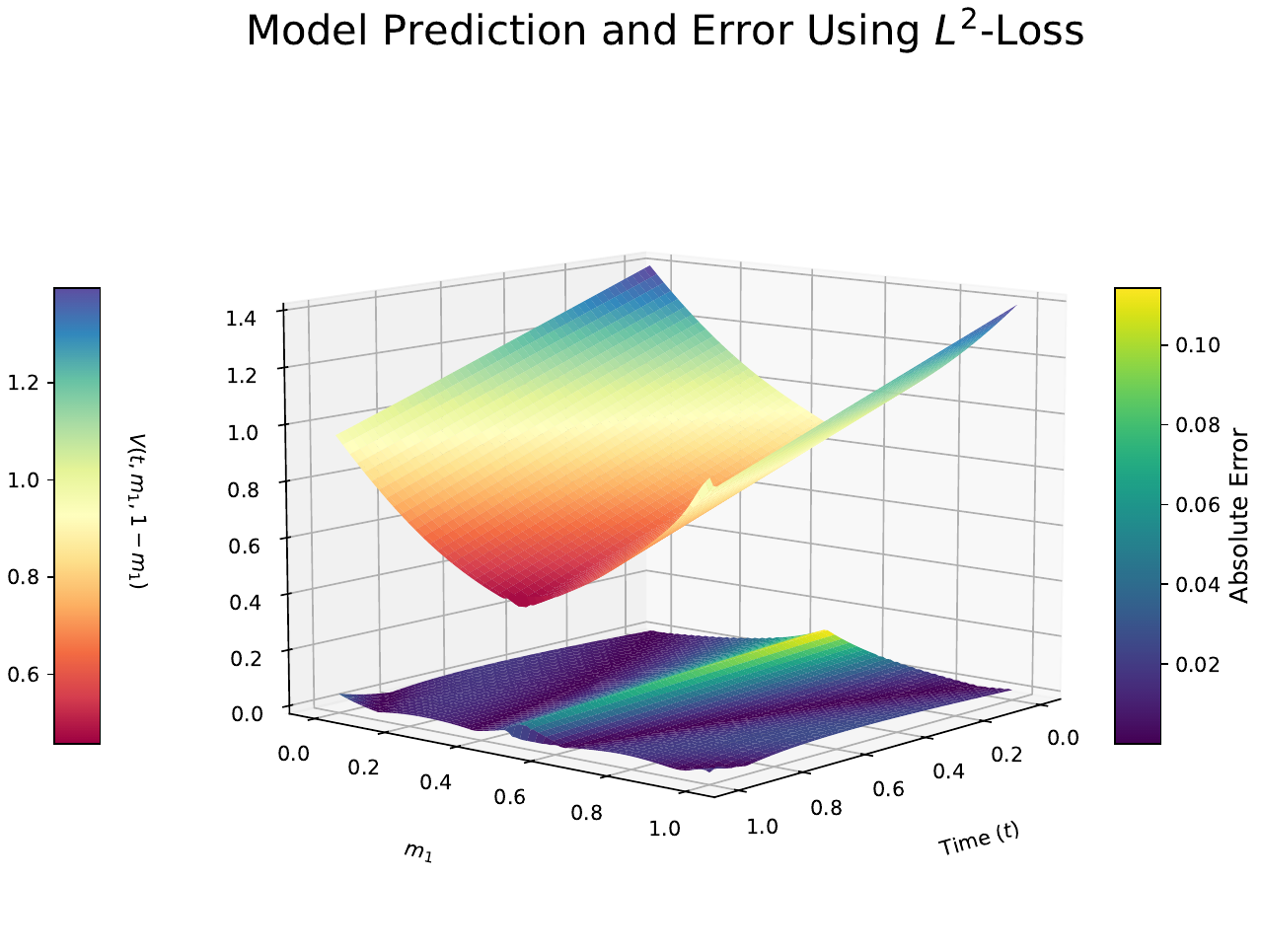}
  \captionof{figure}{\footnotesize Approximate value function, learned using $L^2$-loss, in dimension $d = 2$. The absolute error, measured between the approximate value function and the finite difference solution in Figure~\ref{fig:d=2-true-sol} is displayed on the same plot.}
  \label{fig:d=2-l2surface}
\end{minipage}
\end{figure}

We expect that both metrics are susceptible to local minima using the LSTM-like architecture outlined in \cite{dgm}. Given sufficient training time and appropriate hyperparameter tuning, including the number of samples at each step, the learning rate schedule and the optimizer in use, both algorithms can closely approximate the true terminal condition of the example problem as demonstrated in Figure~\ref{fig:term-condition}. We expect that a deeper neural network, coupled with increased training time, would further improve accuracy. Additionally, it may be possible to obtain a smaller loss by utilizing a quasi-Newton method such as L-BFGS to fine-tune. Recent work suggests that PDE-informed loss functionals suffer from ill-conditioning that can be corrected with such optimization methods \cite{rathore24a}. However, we defer such studies to a future work, having demonstrated the validity of the algorithm with a relatively simple architecture here.

\begin{figure}[htbp]
\centering
\begin{minipage}{.5\textwidth}
  \centering
  \includegraphics[width=\linewidth]{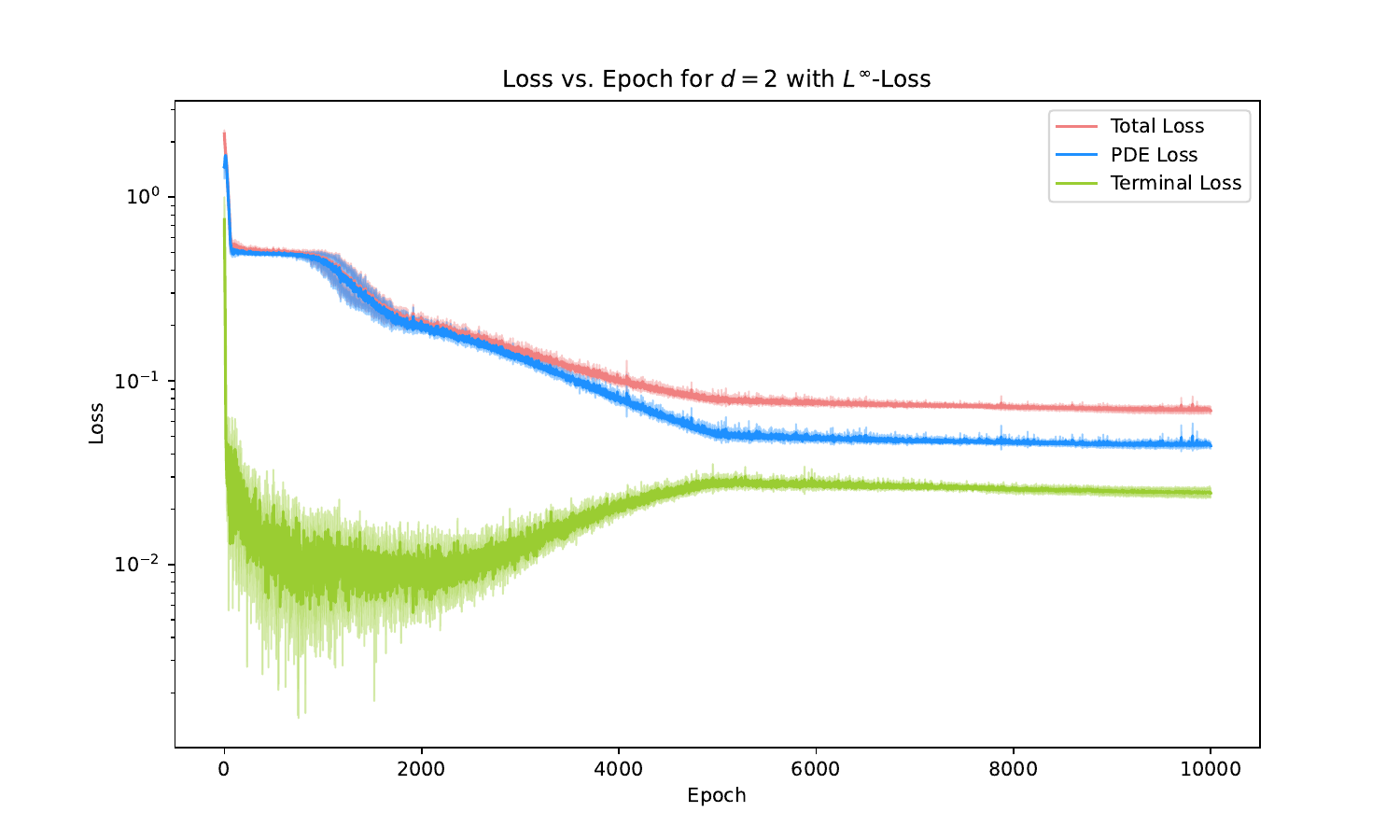}
  \captionof{figure}{\footnotesize $L^\infty$ training loss in dimension $d = 2$. PDE loss refers to the first term in the $L^\infty$-loss from Equation~\eqref{eq:unif-loss}, terminal loss refers to the second term in Equation~\eqref{eq:unif-loss}, and the combined loss represents the entirety of Equation~\eqref{eq:unif-loss}. Shaded regions represent error bars, averaged over five trials.}
  \label{fig:unif-loss}
\end{minipage}%
\begin{minipage}{.5\textwidth}
  \centering
  \includegraphics[width=\linewidth]{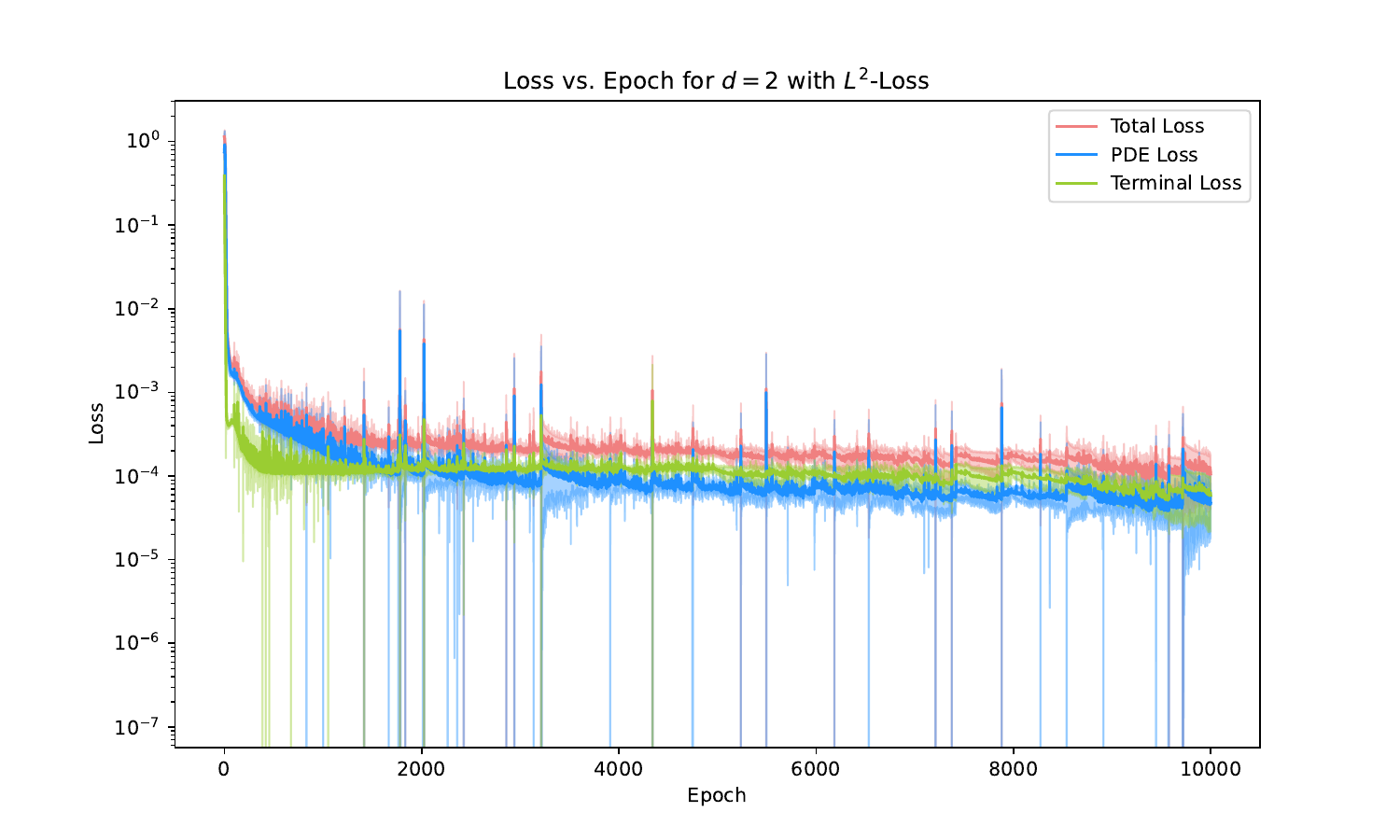}
  \captionof{figure}{\footnotesize $L^2$ training loss in dimension $d = 2$. The total loss, PDE loss, and terminal loss are all defined as in Figure~\ref{fig:unif-loss}. Shaded regions represent error bars, averaged over five trials. Note that the error bars are distorted below the loss curves due to the logarithmic scale of the plot.}
  \label{fig:l2-loss}
\end{minipage}%
\end{figure}

In Figure~\ref{fig:dep-samples-unif}, we demonstrate the relationship between the sample size $M$ and the stability of the loss for the uniform DGM algorithm. Unless otherwise specified, all numerical tests are carried out with $M = 5000$, as the tradeoff between stability and runtime becomes worse as the number of samples exceeds $M = 5000$. As $M$ increases, we observe more stable decreases in loss, as expected. Note that, in Figure~\ref{fig:dep-samples-unif}, smaller values of $M$ seem to yield lower loss values. This phenomenon is potentially misleading due to two reasons. First, with fewer samples $M$ per epoch, points at which the approximation performs poorly are sampled with lower probability, leading to misleadingly low loss values. Second, as one might expect, the loss landscape of the objective appears smoother when taking the maximum over fewer samples. Despite this, training with fewer samples yields a less stable training process that is more susceptible to subpar local minima that do not allow the neural network to accurately approximate the true value function.
\begin{figure}[htbp]
\centering
\begin{minipage}[t]{.49\textwidth}
\vspace{0pt}
\centering
\includegraphics[width=\linewidth]{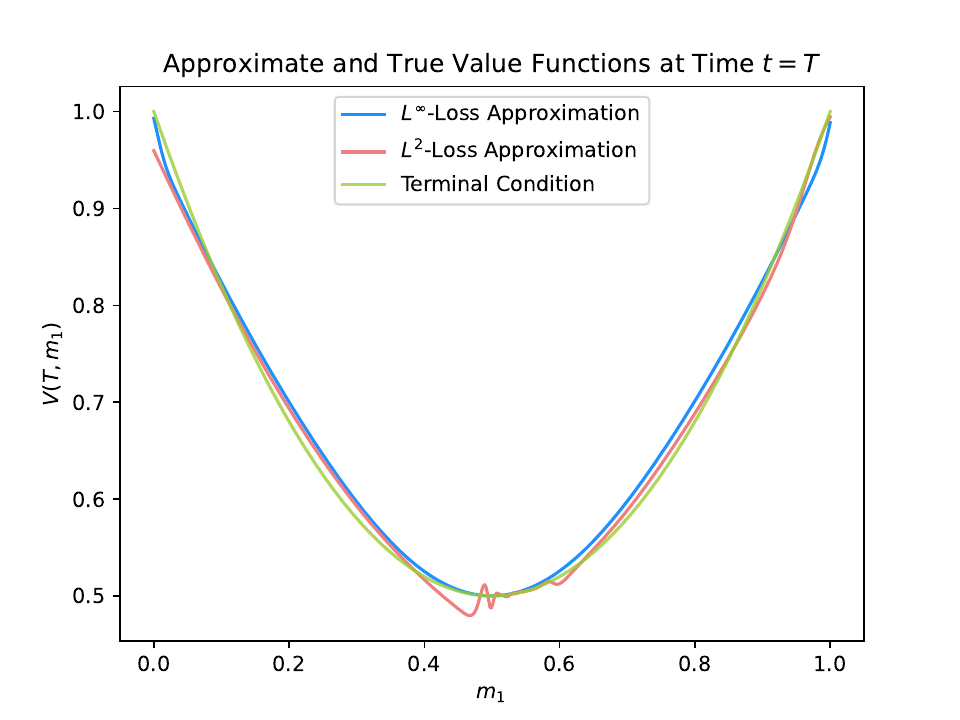}
\captionof{figure}{\footnotesize Approximate and true terminal conditions for both DGM algorithms in dimension $d = 2$.}
\label{fig:term-condition}
\end{minipage}
\begin{minipage}[t]{.49\textwidth}
\vspace{0pt}
\centering
\includegraphics[width=\linewidth]{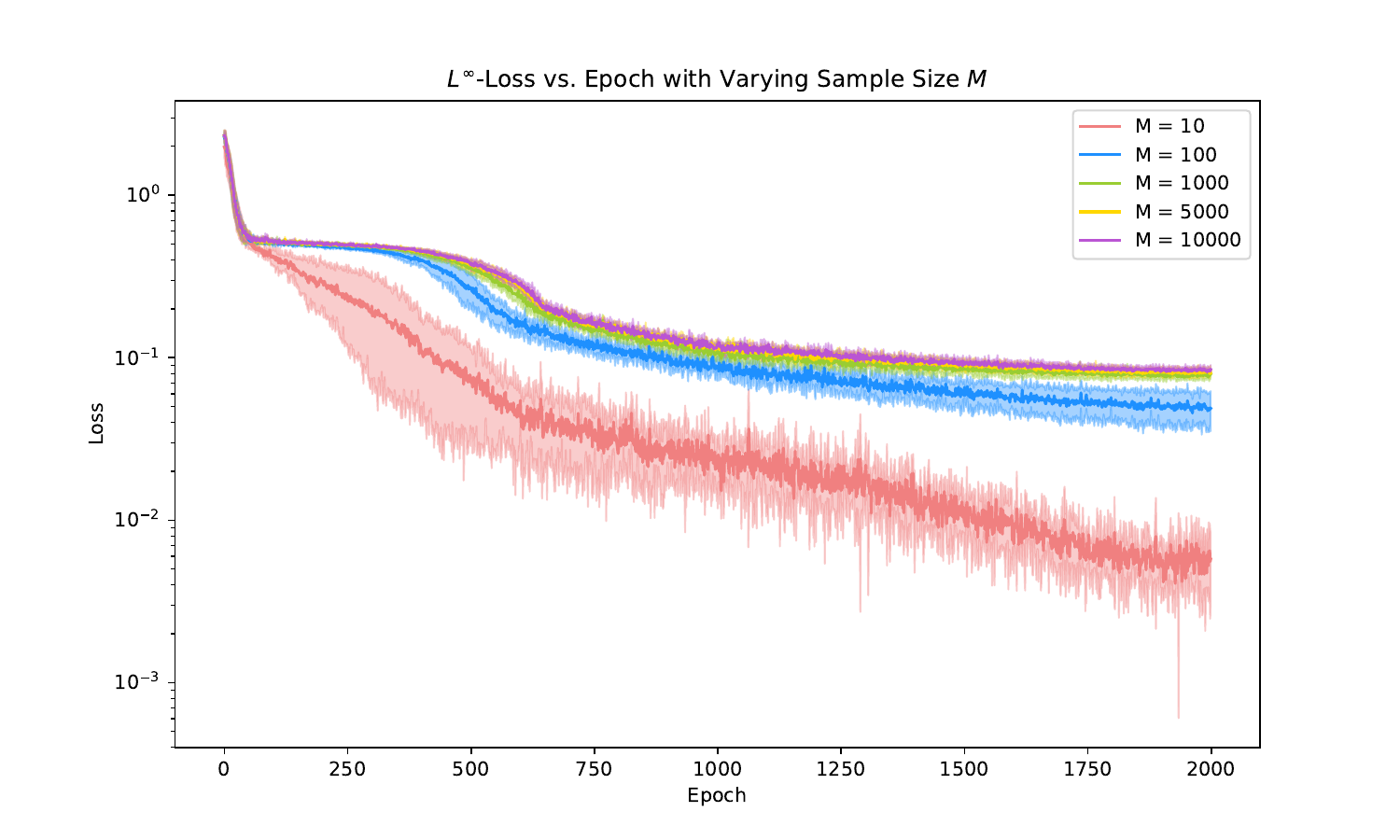}
\captionof{figure}{\footnotesize Dependence of $L^\infty$-loss on number of samples $M$ in dimension $d = 2$. As the number of samples increases, the training process becomes more stable.}
\label{fig:dep-samples-unif}
\end{minipage}
\end{figure}
Finally, in Table \ref{tab:runtimes}, we demonstrate the scalability of the DGM algorithm, implemented in JAX and run with GPU acceleration. By utilizing JAX's built-in auto-differentiation, just-in-time compilation, and GPU acceleration, the DGM algorithm scales well to dimension $d = 200$, with roughly linear increases in runtime past dimension $d = 100$. By runtime, we refer to the time that it takes to train the DGM network with uniform loss for $10000$ epochs, with $10$ gradient descent steps in each epoch, and $M = 5000$ for each epoch.

As the dimension increases, we observe that the loss decreases slightly. This phenomenon could be caused by several factors. Most importantly, as the underlying dimension increases, points at which the neural network poorly approximates the true solution to the HJB equation are sampled with much lower probability. Additionally, while the $L^\infty$-loss functional has both theoretical and practical advantages, a small but nonzero loss still indicates a gap between the approximate solution and the true solution, and the relationship between the loss and the size of this gap may be dimension-dependent. Nonetheless, we expect that our method still approximates the solution to the HJB equation quite well in high dimension, based on its accuracy in dimension $d = 2$ and the broader success of the DGM for solving high-dimensional PDEs, as presented in \cite{dgm}.

\begin{table}[htbp]
\caption{
\small Uniform DGM training times and losses as dimension $d$ increases. PDE loss refers to the first term in the DGM loss from Equation~\eqref{eq:unif-loss}, while terminal loss refers to the second term in Equation~\eqref{eq:unif-loss}. For each dimension $d$, training was repeated five times with five different random seeds, and all reported values are averaged over all runs.}
\label{tab:runtimes}
\vskip 0.15in
\begin{center}
\begin{small}
\begin{sc}
\begin{tabular}{lcccr}
\toprule
Dimension $d$ & Total Training Time ($s$) & Combined Loss & PDE Loss & Terminal Loss \\
\midrule
2 & 315.4 & 0.068 & 0.044 &  0.024 \\
5 & 320.4 & 0.063 & 0.042 & 0.021 \\
10 & 328.2 & 0.042 & 0.025 & 0.016 \\
20 & 335.8 & 0.018 & 0.0096 &  0.0087 \\
50 & 368.8 & 0.0049 & 0.0022 & 0.0027 \\
100 & 448.6 & 0.0019 & 0.00077 & 0.0011 \\ 
200 & 673.2 & 0.00070 & 0.00026 & 0.00044 \\
\bottomrule
\end{tabular}
\end{sc}
\end{small}
\end{center}
\vskip -0.1in
\end{table}
\end{example}
\section{Conclusion}
We present a novel method for solving high-dimensional HJB equations on the simplex arising from MFCPs, complete with a convergence proof our algorithm. Our algorithm, which utilizes an $L^\infty$-loss functional for training rather than a standard $L^2$-loss functional, allows us to leverage the theory of viscosity solutions to prove that if a section of neural network approximators takes the $L^\infty$-loss functional to zero, then the sequence must uniformly converge to the unique viscosity solution of the HJB equation. Notably, our approach does \textit{not} assume boundedness and equicontinuity of the sequence of neural network approximators as in \cite{dgm}, and we thus provide a more general convergence guarantee than that of Sirignano and Spiliopoulos. Additionally, our approach, based on the theory of viscosity solutions, may provide a general framework for neural network-based approaches for solving nonlinear second-order PDEs that: (1) admit a unique $\mathcal{C}^2$ solution, (2) are defined on compact domains, and (3) fall under the class of PDEs that admit viscosity solutions, as defined in \cite{visc-users}. 

Despite the fact that we focus on solving a first-order HJB equation associated with the deterministic MFCP, \cite[Example 1.6]{visc-ishii} shows that the stability properties of viscosity solutions that we leverage in Theorem~\ref{th:unif-value-conv} also hold for second-order HJB equations that arise from stochastic control problems and control of systems of stochastic differential equations. Assuming the existence of a unique $\mathcal{C}^2$ solution, we anticipate that much of our analysis can be extended to show that neural networks can solve such HJB equations on compact domains. An example of such a PDE, to which the DGM has already been successfully applied in practice, can be found in \cite[Section 5]{dgm}. In future work, we aim to relax our regularity assumptions for solutions to the HJB equation and investigate the convergence and performance of deep BSDE, another popular deep learning-based method for solving high-dimensional PDE, when applied to the MFCP.

\vspace{10pt}
\noindent
{\bf Acknowledgements.} 
The authors would like to thank Zihan Guo for his assistance with the numerical experiments presented in this paper. We are also grateful to the Associate Editor Kavita Ramanan from the SIAM Journal on Mathematics of Data Science (SIMODS) and three anonymous referees for their insightful comments, which significantly improved the quality of the manuscript. W. Hofgard is supported by the National Science Foundation Graduate Research Fellowship Program under Grant No. 2146752. This work was also supported by the National Science Foundation (grant DMS-2006305 and DMS-2505998) and by departmental funds from the Department of Mathematics at the University of Michigan, Ann Arbor.

\bibliographystyle{unsrt}
\bibliography{sources}
\normalfont
\end{document}